\documentclass[11pt]{amsart} 
\usepackage{upref,amssymb,eucal} 

\usepackage{diagrams} 
\diagramstyle[small,labelstyle=\scriptstyle,midshaft,PS=dvips] 
\newarrow{Same}===== 
\newarrow{Corresponds}<---> 
\newarrow{Dashonto}....{>>} 
\newarrow{Dashembed}>...> 
 
\theoremstyle{plain} 
\newtheorem{theorem}{Theorem}[section] 
\newtheorem{lemma}[theorem]{Lemma} 
\newtheorem{proposition}[theorem]{Proposition} 
\newtheorem{corollary}[theorem]{Corollary}

\theoremstyle{definition} 
\newtheorem{example}[theorem]{Example} 
 
\newtheorem{definition}[theorem]{Definition} 

\theoremstyle{remark} 
\newtheorem{remark}[theorem]{Remark}

\DeclareMathOperator{\Ext}{Ext} 
\DeclareMathOperator{\shExt}{\mathcal E\mathit{xt}} 
\DeclareMathOperator{\hExt}{\mathbb Ext} 
\DeclareMathOperator{\Hom}{Hom}

\DeclareMathOperator{\coker}{coker} 
 
\DeclareMathOperator{\im}{im} 
 
\DeclareMathOperator{\tr}{tr} 
 
\DeclareMathOperator{\rk}{rk} 
 
\DeclareMathOperator{\Spec}{Spec}

\DeclareMathOperator{\Pic}{Pic}

\DeclareMathOperator{\Pf}{Pf} 
 
\DeclareMathOperator{\gra}{\Gamma}
\newcommand{\graph}[1]{\gra_{#1}}
\DeclareMathOperator{\image}{im}

\newcommand{\rest}[1]{{\mid}_{#1}} 
\newcommand{\sm}[1]{\left(\begin{smallmatrix}#1\end{smallmatrix}\right)} 
\newcommand{\cpx}{{\centerdot}} 
\newcommand{\mult}{{\text{mult}}}

\newcommand{\red}{{\text{\rm red}}} 
\newcommand{\linsys}[1]{\lvert #1 \rvert} 
\newcommand{\cond}[1]{\ref{cond.#1}} 
\newcommand{\hyp}{q_h}
 
\renewcommand\O{{\mathcal O}} 
\renewcommand\P{{\mathbb P}}

\begin{document}

\title[Lagrangian Subbundles and Subcanonical Subschemes] 
{Lagrangian Subbundles and Codimension $3$ Subcanonical Subschemes} 
\author{David Eisenbud}  
\address{Department of Mathematics\\ University of California, 
Berkeley\\ Berkeley CA 94720} 
\email{de@msri.org} 
\urladdr{http://www.msri.org/people/staff/de/}

\author{Sorin Popescu}
\address{Department of Mathematics\\
Columbia University\\
New York, NY  10027}
\email{psorin@math.columbia.edu}
\urladdr{http://www.math.columbia.edu/\~{}psorin/}
\author{Charles Walter} 
\address{Laboratoire J.-A.\ Dieudonn\'e (UMR CNRS 6621)\\ 
Universit\'e de Nice -- Sophia Antipolis\\ 06108 Nice Cedex 02\\ 
France} 
\email{walter@math.unice.fr} 
\urladdr{http://math1.unice.fr/\~{}walter/} 
\thanks{Partial support for the authors during the preparation of 
this work was provided by the NSF. The authors are also 
grateful to MSRI Berkeley and the University of Nice Sophia-Antipolis 
for their hospitality.}


\begin{abstract} 
  We show that a Gorenstein subcanonical codimension $3$ subscheme
  $Z\subset X=\P^N$, $N\ge 4$, can be realized as the locus along
  which two Lagrangian subbundles of a twisted orthogonal bundle meet
  degenerately, and conversely. We extend this result to singular $Z$
  and all quasiprojective ambient schemes $X$ under the necessary
  hypothesis that $Z$ is {\it strongly subcanonical\/} in a sense
  defined below. A central point is that a pair of Lagrangian
  subbundles can be transformed locally into an alternating map. In
  the local case our structure theorem reduces to that of
  Buchsbaum-Eisenbud \cite{BE} and says that $Z$ is Pfaffian.
   
  We also prove codimension one symmetric and skew-symmetric analogues
  of our structure theorems.
\end{abstract} 
 
\maketitle

\label{introduction} 
 
Smooth subvarieties of small codimension $Z \subset X=\mathbb P^N$
have been extensively studied in recent years, especially in relation
to Hartshorne's conjecture that a smooth subvariety of sufficiently
small codimension in $\mathbb P^N$ is a complete intersection.
Although the conjecture remains open, any smooth subvariety $Z$ of
small codimension in $\mathbb P^N$ is known, by a theorem of Barth,
Larsen, and Lefschetz, to have the weaker property that it is {\em
subcanonical\/} in the sense that its canonical class is a multiple of
its hyperplane class.

More generally, a subscheme $Z$ of a nonsingular Noetherian scheme $X$
is said to be subcanonical if $Z$ is Gorenstein and its canonical
bundle is the restriction of a bundle on $X$. There is a natural
generalization to an arbitrary (possibly singular) scheme $X$
(see below).

In this paper we give a structure theorem for subcanonical 
subschemes of codimension $3$ in $\P^N$ and generalize it to
subcanonical subschemes of codimension $3$ in an
arbitrary quasiprojective scheme $X$ satisfying a mild extra
cohomological condition (strongly subcanonical subschemes). 
The construction works even without the quasiprojective
hypothesis.
 
There are well known theorems describing the local structure of
Gorenstein subschemes of nonsingular Noetherian schemes in
codimensions $\leq 3$.  In codimensions $1$ and $2$ all Gorenstein
subschemes are locally complete intersections.  These results have
been globalized: If $X$ is nonsingular, any $Z \subset X$ of
codimension $1$ is the zero locus of a section of a line bundle, while
a subcanonical $Z \subset X$ of codimension $2$ is the zero locus of a
section of a rank $2$ vector bundle if a certain obstruction in
cohomology vanishes (as explained below).  In both cases $\mathcal
O_Z$ has a symmetric resolution by locally free $\mathcal
O_X$-modules.
 
In codimension $3$ both the local and the global cases become more 
complicated.  Locally, a Gorenstein subscheme of codimension $3$ need 
not be a locally complete intersection.  Rather, Buchsbaum and Eisenbud 
\cite{BE} showed that such a subscheme is cut out locally by the 
submaximal Pfaffians of an alternating matrix appearing in a minimal 
free resolution.  Okonek \cite{Okonek} asked whether this local result 
could be generalized to show that codimension $3$ subcanonical schemes 
are cut out by the Pfaffians of an alternating map of vector bundles. 
Walter \cite{Walter} gave a positive answer to Okonek's question in 
$\mathbb P^n$ under a mild additional hypothesis, but left open the 
question of whether this hypothesis is always satisfied. 
 
In our paper \cite{EPW2} we will show that {\em not every 
subcanonical subscheme of codimension $3$ in $\mathbb P^n$ is 
Pfaffian,\/} settling Okonek's question negatively.  But in the 
present paper we show that a different way of looking at the Pfaffian 
construction does generalize, and gives the desired structure theorem 
for all subcanonical subschemes of codimension $3$. (The question as 
to which subschemes are Pfaffian can be answered in the derived Witt 
group of Balmer \cite{Balmer}; see Walter \cite{Walter2}.) 
 
In this paper
a closed subscheme $Z \subset X$ of a Noetherian scheme is called
{\em subcanonical of codimension $d$} if it satisfies two 
conditions:  
\begin{enumerate} 
\renewcommand\theenumi{\Alph{enumi}} 
\renewcommand\labelenumi{\theenumi)} 
\setlength{\itemsep}{0.6ex plus 0.6ex minus 0ex} 
\item \label{cond.1} The subscheme $Z$ is {\em relatively 
Cohen-Macaulay of codimension $d$} in $X$, i.e.\ $\shExt^i_{\mathcal 
O_X} (\mathcal O_Z, \mathcal O_X) = 0$ for all $i \neq d$, and 
\item \label{cond.2} There exists a line bundle $L$ on $X$ such that 
the relative canonical sheaf $\omega_{Z/X}:= \shExt^d_{\mathcal 
O_X}(\mathcal O_Z,\O_X)$ is isomorphic to the restriction of $L^{-1}$ 
to $Z$. 
\end{enumerate} 
These conditions are not enough for the Serre correspondence in
codimension $2$, nor for our structure theorem in codimension $3$.
Condition \cond2 asserts the existence of an isomorphism
\[
\eta : \mathcal O_Z \xrightarrow{\sim} \omega_{Z/X}(L) = 
\shExt^d_{\mathcal O_X}(\mathcal O_Z,L) 
\]
which one can think of as an $\eta \in H^0(\shExt^d_{\mathcal O_X} 
(\mathcal O_Z,L))= \Ext^d_{\mathcal O_X} (\mathcal O_Z,L)$.  In the 
Yoneda Ext, this $\eta$ defines a class of ``resolutions of $\mathcal 
O_Z$ by coherent sheaves'' 
\[ 
0 \to L \to \mathcal F_{d-1} \to \dotsb \to \mathcal F_1 \to \mathcal 
F_0 \to \mathcal O_Z \to 0. 
\] 
In our structure theorem we require such resolutions with
$\mathcal F_0, \dots, \mathcal F_{d-1}$ 
locally free.  Thus we will need the condition
\begin{enumerate} 
\renewcommand\theenumi{\Alph{enumi}} 
\renewcommand\labelenumi{\theenumi)} 
\setlength{\itemsep}{0.6ex plus 0.6ex minus 0ex} 
\addtocounter{enumi}{2} 
\item \label{cond.3} The $\mathcal O_X$-module $\mathcal O_Z$ is of 
  finite local projective dimension (necessarily equal to the 
  codimension $d$).
\end{enumerate}
This condition holds automatically if the ambient scheme $X$
is nonsingular.

We also need $\mathcal F_0 = \mathcal O_X$, which means that we want
$\eta \in \Ext^d_{\mathcal O_X}(\mathcal O_Z,L)$ to lift to
$\Ext^{d-1}_{\mathcal O_X}(\mathcal I_Z,L)$.  Since these two groups
are joined by a map in the long exact sequence obtained by applying
$\Ext^*_{\mathcal O_X}({-},L)$ to the short exact sequence
$0\to\mathcal I_Z\to \mathcal O_X\to\mathcal O_Z\to 0$, we see that
the lifting exists if and only if $\eta\in\Ext^d_{\mathcal O_X}
(\mathcal O_Z,L)$ goes to $0$ in $\Ext^d_{\mathcal O_X} (\mathcal
O_X,L)\cong H^d(X,L)$.  We are thus led to the condition:
\begin{enumerate} 
\renewcommand\theenumi{\Alph{enumi}} 
\renewcommand\labelenumi{\theenumi)} 
\setlength{\itemsep}{0.6ex plus 0.6ex minus 0ex} 
\addtocounter{enumi}{3} 
\item \label{cond.4} The isomorphism class $\eta \in \Ext^d_{\mathcal 
O_X} (\mathcal O_Z,L)$ of \eqref{eta} goes to zero under the map 
\[
\Ext^d_{\mathcal O_X}(\mathcal O_Z,L) \to \Ext^d_{\mathcal 
O_X}(\mathcal O_X,L) = H^d(X,L) 
\]
induced by the surjection $\mathcal O_X \to \mathcal O_Z$.  
\end{enumerate} 
Condition \cond4 holds automatically if $H^d(X,L) = 0$. 
This is the case if $X = \mathbb P^n$ with $n \geq d+1$, or if $X$ is 
an affine scheme.  In addition, if the ambient scheme $X$ is a 
Gorenstein variety over a field $k$, then condition \cond4 can be put 
into a dual form which looks more natural.  For in that case $Z 
\subset X$ is subcanonical (of dimension $r$) if and only if it is 
Cohen-Macaulay and there exists a line bundle $M$ on $X$ such that 
$\omega_Z \cong M \rest Z$, and condition \cond4 holds if and only if 
the following composite map vanishes 
\begin{equation} 
\label{cond.dual} 
H^r(X,M) \xrightarrow{\text{\rm rest}} H^r(Z,M \rest Z) \xrightarrow 
[\cong] {\eta} H^r(Z,\omega_Z) \xrightarrow{\tr} k. 
\end{equation} 
In these terms we may give the central definition of this paper:

\begin{definition}
\label{strong.def}
A subscheme $Z\subset X$ is
{\em strongly subcanonical} if it satisfies 
conditions \cond1-\cond4.
\end{definition}

The Serre construction shows that a subscheme of codimension $2$ is
the zero locus of a rank $2$ vector bundle if and only if it is
strongly subcanonical (Griffiths-Harris \cite{GH} Proposition 1.33,
Vogelaar \cite{Vo} Theorem 2.1, and B\u anic\u a-Putinar \cite{BP} \S
2.1 state variants of condition \cond4 explicitly).

Our main results show that a codimension $3$ subscheme $Z$ of a
quasiprojective scheme $X$ is strongly subcanonical if and only if it
can be expressed as an appropriate ``Lagrangian degeneracy scheme,''
defined as follows: Let $\mathcal V$ be a vector bundle on $X$ of even
rank $2n$ equipped with a nonsingular quadratic form $q$ with values
in a line bundle $L$.  Let $\mathcal E$ and $\mathcal F$ be a pair of
Lagrangian subbundles of $(\mathcal V,q)$ (i.e.\ totally isotropic
subbundles of rank $n$).  It is then well known that $\dim [ \mathcal
E(x) \cap \mathcal F(x) ]$ is locally constant modulo $2$.
 
Now suppose that $m$ is an integer such that $\dim [ \mathcal E(x) 
\cap \mathcal F(x) ] \equiv m \pmod 2$ for all $x \in X$.  Then there 
is a degeneracy locus which as a set is given by 
\[ 
Z_m(\mathcal E,\mathcal F)_\red := \left \{ x\in X \mid \dim_{k(x)} 
  [\mathcal E(x) \cap \mathcal F(x)] \geq m \right \}. 
\] 
In \S\ref{sect.degen} we will define a scheme structure on 
this set in roughly the following manner.  Using the data $\mathcal E, 
\mathcal F \subset (\mathcal V,q)$ one defines a composite map 
\[ 
\lambda:\quad \mathcal E\rightarrow \mathcal V\cong \mathcal V^*(L) 
\rightarrow \mathcal F^*(L) 
\] 
such that $\ker(\lambda (x)) = \mathcal E(x) \cap \mathcal F(x)$ for 
all $x \in X$.  
Even if $\mathcal E\cong \mathcal F$ the map $\lambda$ may not be
alternating, but (perhaps after modifying $\lambda$ slightly to make
it have even rank everywhere) we will show that it is possible to find
local bases in which the matrix of $\lambda$ is alternating
(Proposition \ref{prop.loc.alt}).  Although these alternating
matrices do not glue together, they are sufficiently compatible that
we can define $Z_m(\mathcal E,\mathcal F)$ as the locus defined by
their Pfaffians of order $\rk(\mathcal E)-m+2$.  This scheme structure
is natural in the sense that in a suitably generic setting it is
reduced, and it is stable under base change.

The following structure theorem for strongly subcanonical codimension
$3$ subschemes collects the main results of this paper:
\begin{theorem}
\label{main.result}
Let $X$ be a quasiprojective scheme over a Noetherian ring,
and let $Z\subset X$ be a closed subscheme of grade $3$. 
The following conditions are equivalent\textup:

\textup{(a)} $Z\subset X$ is strongly subcanonical.

\textup{(b)} There exists a twisted orthogonal bundle $(\mathcal V, q)$
and Lagrangian subbundles $\mathcal E,\mathcal F\subset (\mathcal V,q)$,
with  $\dim_{k(x)} \bigl[ \mathcal E(x) 
\cap \mathcal F(x) \bigr]$ odd for all $x 
\in X$, such that $Z=Z_3(\mathcal E, \mathcal F)$.

\textup{(c)} There exists a vector bundle $\mathcal F$, a line bundle
$L$, and a Lagrangian subbundle $\mathcal E$ of the hyperbolic bundle
$\mathcal F\oplus\mathcal F^*(L)$, such that the composite map
\[
\lambda: \mathcal E \hookrightarrow \mathcal F\oplus \mathcal F^*(L)
\twoheadrightarrow \mathcal F^*(L)
\]
has kernel of odd rank and such that $Z=Z_3(\mathcal E, \mathcal F^*(L))$.

\textup{(d)} $Z$ has symmetrically quasi-isomorphic locally
free resolutions 
\[
\begin{diagram}
0 & \rTo & L & \rTo & \mathcal H & \rTo ^\psi & \mathcal G & \rTo & 
\mathcal O_X & \rTo & \mathcal O_Z & \rTo & 0 \\ 
&& \dSame && \dTo <\phi && \dTo >{\phi^*} && \dSame && \dTo >\eta 
<\cong \\ 
0 & \rTo & L & \rTo & \mathcal G^*(L) & \rTo ^{-\psi^*} & \mathcal 
H^*(L)& \rTo & \mathcal O_X & \rTo & \shExt^3_{\mathcal O_X} (\mathcal 
O_Z, L) & \rTo & 0 
\end{diagram} 
\]
where $L$ is a line bundle on $X$, and $\phi^*\psi : \mathcal E \to
\mathcal E^*(L)$ is an alternating map.
\end{theorem}

The structure theorem will be proved in several parts: see Theorems
\ref{second.v}, \ref{first.v}, and \ref{converse}.

One way to look at the structure theorem is as follows.  The existence
of the symmetric isomorphism
\begin{equation}
\label{eta} 
\eta : \mathcal O_Z \xrightarrow{\sim} \shExt^3_{\mathcal
O_X}(\mathcal O_Z,L)
\end{equation} 
means that there should be a symmetric isomorphism in the derived
category from the locally free resolution of $\mathcal O_Z$
\begin{equation} 
\label{first.resol} 
0 \to L \xrightarrow{\,q\,} \mathcal E \xrightarrow{\,\psi\,} \mathcal 
G \xrightarrow{\,p\,} \mathcal O_X \to \mathcal O_Z \to 0
\end{equation} 
into its twisted shifted dual.  In general, morphisms in the derived
category are complicated objects involving homotopy classes of maps
and a calculus of fractions.  Nevertheless, in Theorem \ref{converse}
we show that there exist such locally free resolutions of $\mathcal
O_Z$---which depend on the choice of $\eta$---for which the symmetric
isomorphism in the derived category is induced by a symmetric
chain map which is a quasi-isomorphism 
and therefore becomes an
isomorphism in the derived category.  Okonek's Pfaffian subschemes
correspond to situations where this quasi-isomorphism is an
isomorphism.
 
The philosophy that (skew)-symmetric sheaves should have locally free
resolutions that are (skew)-symmetric up to quasi-isomorphism is also
pursued in our paper \cite{EPW3} and in Walter \cite{Walter2}.  The
former deals primarily with methods for constructing explicit locally
free resolutions for (skew)-symmetric sheaves on $\mathbb P^n$.  The
latter studies the obstructions (in Balmer's derived Witt groups
\cite{Balmer}) to the existence of a genuinely (skew)-symmetric
resolution.
 
The results of this paper give a full characterization of codimension
$3$ subcanonical subschemes.  In \cite{EPW2} we use this machinery to
construct various geometric examples of subcanonical
subschemes of codimension $3$ which are not Pfaffian.
 
Porteous-type formulas for the fundamental classes of degeneracy loci
for skew-symmetric maps $\phi : \mathcal E \to \mathcal E^*(L)$ were
found by Harris-Tu \cite{HarrisTu}, J\'ozefiak-Lascoux-Pragacz
\cite{JLP}, and Pragacz \cite{Pragacz}.  Harris asked for similar
formulas for degeneracy loci related to pairs of Lagrangian
subbundles, and they were provided by Fulton \cite{FultonHirz}
\cite{FultonJDG} and Pragacz-Ratajski \cite{PR} (see Fulton-Pragacz
\cite{FP} for more details).  (A scheme structure on these degeneracy
loci and their generalizations with isotropic flag conditions can be
defined in a manner similar to \eqref{scheme.str} below.)
 
Fulton and Pragacz (\cite{FP} \S 9.4) also ask whether one can find 
``natural'' resolutions for the structure sheaves of these kinds of 
symmetric and skew-symmetric degeneracy loci.  From such a resolution 
one can read off formulas in $K_0(X)$.  Theorem \ref{second.v} 
provides an explicit answer in one simple case. 
 
\subsection{Structure of the paper} 

In Sections \ref{sect.lagr} and \ref{sect.degen} we review basic facts
about Lagrangian subbundles of twisted orthogonal bundles and define
the scheme structure on the degeneracy loci $Z_m(\mathcal E,\mathcal
F)$.  In Section \ref{Subcanonical} we prove that Lagrangian
degeneracy loci of codimension $3$ are strongly subcanonical (Theorem
\ref{second.v}).  In Section \ref{lagrange} we discuss ``split''
Lagrangian degeneracy loci (Theorem \ref{first.v}) which are often
more practical for constructing codimension $3$ subcanonical
subschemes.  The computation of local equations for these degeneracy
loci is discussed in Section \ref{local.eq}.
 
In Section \ref{converse.theorem} we complete the proof of the
structure Theorem \ref{main.result} by showing that strongly
subcanonical subschemes of codimension $3$ are split Lagrangian
degeneracy loci (Theorem \ref{converse}).  Sections \ref{sect.points}
and \ref{examples} discuss at length various examples of codimension
$3$ subcanonical subschemes, particularly the case of points in
$\mathbb P^3$.  Further examples can be found in our paper
\cite{EPW2}.
 
Finally, in Section \ref{codimension.one} we prove codimension one
symmetric and skew-symmetric analogues of all previous results.  In
particular we state Casnati-Catanese's structural result (\cite{CC}
Remark 2.2) and give an example of a self-linked threefold of degree
$18$ in $\mathbb P^5$ which does not have a symmetric resolution
because the parity condition fails.

\subsection{Acknowledgements} 
 
We are grateful for many useful discussions to Fabrizio Catanese,
Igor Dolgachev, Bill Fulton, Joe Harris, Andr\'e Hirschowitz,
and Roberto Pignatelli. 
 
The second and third authors would like to thank the Mathematical 
Sciences Research Institute in Berkeley for its support while part of 
this paper was being written.  The first and second authors also thank 
the University of Nice -- Sophia Antipolis for its hospitality.   
 
Many of the diagrams were set using the {\tt diagrams.tex} 
package of Paul Taylor.

\section{Quadratic forms on vector bundles} 
\label{sect.lagr} 
 
In this section we recall the basic definitions of twisted orthogonal 
bundles and of Lagrangian subbundles.  The definitions and results can 
be found in many standard references such as Fulton-Pragacz \cite{FP} 
Chap.\ 6, Knus \cite{Knus}, and Mukai \cite{Mukai} \S 1.

\subsection{Quadratic forms} 
 
Suppose that $V$ is a finite-dimensional vector space over a field 
$k$.  (We impose no restrictions on $k$; it
may have characteristic 2, and need not be algebraically 
closed.)  A quadratic form on $V$ is a 
homogeneous quadratic polynomial in the linear forms on $V$, i.e.\ a 
member $q \in S^2(V^*)$.  The symmetric bilinear form $b: V\times V 
\to k$ associated to $q$ is given by the formula 
\begin{equation} 
\label{ass.bilin.form} 
b(x,y) := q(x+y) - q(x) - q(y). 
\end{equation} 
The quadratic form $q$ is {\em nondegenerate} if $b$ is a perfect 
pairing.   
 
Now suppose that $\mathcal V$ is a locally free sheaf of constant 
finite rank over a scheme $X$.  A quadratic form on $\mathcal V$ with 
values in a line bundle $L$ is a global section $q$ of $S^2(\mathcal 
V^*) \otimes L$.  Such a quadratic form is {\em nonsingular} if the 
induced symmetric bilinear form is a perfect pairing.  Equivalently a 
quadratic form $q$ on $\mathcal V$ is nonsingular if for each point $x 
\in X$ the induced quadratic form $q(x)$ on the fiber vector space 
$\mathcal V(x)$ is nondegenerate.  A {\em twisted orthogonal bundle} 
on $X$ is a vector bundle $\mathcal V$ equipped with a nonsingular 
quadratic form $q$ with values in some line bundle $L$.

\subsection{Lagrangian subbundles}  
 
If $V$ is a vector space of even dimension $2n$ equipped with a 
nondegenerate quadratic form, then a {\em Lagrangian subspace} $E 
\subset (V,q)$ is a subspace of $V$ of dimension $n$ such that $q 
\rest{E} \equiv 0$.  If the characteristic is $\neq 2$, then $E 
\subset (V,q)$ is Lagrangian if and only if $E = E^\perp:= \{ x\in V 
\mid b(x,y)=0 \text{ for all $y \in E$} \}$.  But in characteristic 
$2$ this condition is necessary but not sufficient for $q$ to vanish 
on $E$, i.e.\ for $E$ to be Lagrangian. 
 
Similarly, a {\em Lagrangian subbundle} $\mathcal E \subset (\mathcal 
V,q)$ of a twisted orthogonal bundle of even rank $2n$ is a subbundle 
(with locally free quotient sheaf) of rank $n$ such that $q 
\rest{\mathcal E} \equiv 0$.

The following result is well known (cf.\ Bourbaki \cite{Bourbaki} \S 6 
ex.\ 18(d), Mumford \cite{Mumford}, Mukai \cite{Mukai} Proposition 
1.6). 
 
\begin{proposition} 
\label{constant.parity} 
If $\mathcal E$ and $\mathcal F$ are Lagrangian subbundles of a 
twisted orthogonal bundle over a scheme $X$, then the function on $X$ 
given by $x \mapsto \dim_{k(x)} \bigl[ \mathcal E(x) \cap \mathcal 
F(x) \bigr]$ is locally constant modulo $2$. 
\end{proposition}

\subsection{Hyperbolic bundles} 
 
If $\mathcal F$ is any vector bundle of constant rank, and $L$ is any
line bundle, then then $\mathcal F\oplus \mathcal F^*(L)$
may be endowed with the {\em hyperbolic quadratic form} $\hyp(e \oplus
\alpha) := \alpha(e)$ with values in $L$.  (This $\hyp$ is bilinear on
$\mathcal F \times \mathcal F^*(L)$ but quadratic on $\mathcal F
\oplus \mathcal F^*(L)$.)  The associated hyperbolic symmetric
bilinear form has matrix $\sm{0&I\\I&0}$.

We will use the following notation for graph subbundles.  If $\psi :
\mathcal A \to \mathcal B$ and $\alpha : \mathcal B \to \mathcal A$
are morphisms of vector bundles, then we write
\begin{align*}
\graph{\psi} & := \image \left(\mathcal A \rInto^{\sm{1 \\ \psi}}
\mathcal A \oplus \mathcal B \right), &
\graph{\alpha} & := \image \left(\mathcal B \rInto^{\sm{\alpha \\ 1}}
\mathcal A \oplus \mathcal B \right).
\end{align*}
These graphs are to be regarded as subbundles of $\mathcal A \oplus
\mathcal B$.

\begin{lemma} 
\label{graph.alt} 
A subbundle $\mathcal E \subset (\mathcal F \oplus \mathcal F
^*(L),q_h)$ is a Lagrangian subbundle complementary to the direct
summand $\mathcal F^*(L)$ if and only if there is an alternating map
$\zeta: \mathcal F \to \mathcal F^*(L)$ such that $\mathcal E =
\graph{\zeta}$.
\end{lemma}

Any Lagrangian subbundle of a twisted orthogonal bundle over an
\textbf{affine} scheme has a Lagrangian complement (cf.\ \cite{Knus}
Remark I.5.5.4), although this is not always true over a general
scheme.  However, if a Lagrangian subbundle $\mathcal F \subset
(\mathcal V,q)$ has a Lagrangian complement $\mathcal M$, then the
symmetric bilinear form induces a natural isomorphism $\mathcal M
\cong \mathcal F^*(L)$.  This defines an isometry
\begin{equation}
\label{isom.hyp}
(\mathcal V, q) \rTo ^{\phi_{\mathcal F,\mathcal M}} _\cong
(\mathcal F \oplus \mathcal F^*(L),\hyp)
\end{equation}
which is the identity on $\mathcal F$ and which identifies the
complementary Lagrangian subbundles $\mathcal F, \mathcal M \subset
\mathcal V$ with the two direct summands of $\mathcal F \oplus
\mathcal F^*(L)$.  The previous lemma has the following corollary.

\begin{corollary}
\label{common.compl}
If $\mathcal E, \mathcal F \subset (\mathcal V,q)$ are Lagrangian
subbundles with a common Lagrangian complement $\mathcal M$, then
there is an alternating map $\zeta : \mathcal F \to \mathcal F^*(L)$
such that $\phi_{\mathcal F, \mathcal M} (\mathcal E) = \Gamma_\zeta$.
\end{corollary}

\section{Locally alternating maps and Lagrangian degeneracy loci} 
\label{sect.degen} 

In this section we show how to use Corollary \ref{common.compl} to
define scheme-theoretic degeneracy loci for pairs of Lagrangian
subbundles of a twisted orthogonal bundle which generalize the
degeneracy loci for alternating maps defined by ideals of Pfaffians.
We also show how to turn a pair of Lagrangian subbundles into a
locally alternating map.  Several steps are required, in order to make
sure that common Lagrangian complements exist locally, and to show
that our degeneracy loci are independent of the choice of common
Lagrangian complement. Our scheme structure defined
by local equations coincides with that given by a universal
construction in De Concini-Pragacz \cite{DeConciniPragacz}.

\subsection{Existence of local common Lagrangian complements}
The result we need is is standard if the residue field
is infinite, but if the residue field is very small, care is
required. 
We are interested in when two Lagrangian subbundles $\mathcal E,
\mathcal F$ of a twisted orthogonal bundle $(\mathcal V,q)$ have a
common Lagrangian complement locally.  If one recalls that an
even-dimensional quadratic vector space $(\mathcal V(x), q(x))$ has
two families of Lagrangian subspaces, and that in order for two
Lagrangian subspaces to have a common Lagrangian complement they must
lie in the same family, and that this is measured by the dimension of
the intersection modulo $2$, we see that in order for $\mathcal E$ and
$\mathcal F$ to have a common Lagrangian complement to $\mathcal E$
and $\mathcal F$, we must have $\dim_{k(x)} \left[ \mathcal E(x) \cap
\mathcal F(x) \right] \equiv \rk (\mathcal E) \pmod 2$.  We now show
that locally this condition is also sufficient.

\begin{proposition} 
\label{prop.com.comp} 
Let $\mathcal E, \mathcal F \subset (\mathcal V,q)$ be Lagrangian
subbundles of a twisted orthogonal bundle on a scheme $X$.  Suppose
that $\dim_{k(x)} \bigl[ \mathcal E(x) \cap \mathcal F(x) \bigr]
\equiv \rk(\mathcal E) \pmod 2$ for all $x \in X$.  Then any 
$x \in X$ has a neighborhood $U$ over which $\mathcal E \rest U$ and
$\mathcal F \rest U$ have a common Lagrangian complement $\mathcal
M_U$.
\end{proposition} 

\begin{proof}
It is easy to see that any common Lagrangian complement at $x$ extends
to a common Lagrangian complement in a neighborhood $U$.  The
existence at $x$ of such a complement is standard if the residue field
is infinite.  The following lemma deals with the remaining case:
\renewcommand{\qed}{}
\end{proof}

\begin{lemma} 
\label{three.lagr} 
Suppose that $q$ is a nondegenerate quadratic form on an 
even-dimensional vector space $V$, and that $U, U' \subset (V,q)$ are 
two Lagrangian subspaces such that $\dim (U \cap U') \equiv \dim (U) 
\pmod 2$.  Then there exists a Lagrangian subspace $L \subset (V,q)$ 
complementary to $U$ and to $U'$. 
\end{lemma} 
 
\begin{proof} 
Let $K = U \cap U'$.  Then $U = U^\perp \subset K^\perp$, and 
similarly $U' \subset K^\perp$.  On dimensional grounds, we must 
indeed have $U + U' = K^\perp$.  As a result $U/K$ and $U'/K$ are 
complementary Lagrangian subspaces of $K^\perp/K$.  Moreover, by 
hypothesis they are even-dimensional. 
 
Let $f_1,\dots,f_{2m}$ be a system of vectors in $U$ mapping onto a 
basis of $U/K$.  Since $U/K$ and $U'/K$ are complementary Lagrangian 
subspaces of $K^\perp/K$, the symmetric bilinear form $b$ associated 
to $q$ induces a perfect pairing between them.  So there exists a 
system of vectors $g_1, \dots, g_{2m}$ in $U'$ such that $b(f_i,g_j) = 
\delta_{ij}$ for all $i,j$.   
 
Let $N$ be the subspace spanned by the $f_i$ and $g_j$.  Then $q 
\rest{N}$ is nondegenerate, so there is an orthogonal direct sum 
decomposition $V = N \oplus N^\perp$ such that $q\rest{N}$ and 
$q\rest{N^\perp}$ are both nondegenerate.  Moreover, $K \subset 
(N^\perp, q \rest{N^\perp})$ is a Lagrangian subspace, for which there 
exists a complementary Lagrangian subspace $P$ by our previous 
remarks.  Let $p_1, \dots, p_r$ be a basis of $P$.  One may now check 
that 
\[ 
f_1+g_2, f_2-g_1, f_3+g_4, f_4-g_3, \dots, f_{2m-1}+g_{2m}, 
f_{2m}-g_{2m-1}, p_1, \dots, p_r 
\] 
form a basis for a Lagrangian subspace $L \subset (V,q)$ complementary 
to both $U$ and $U'$.  
\end{proof}

The following example shows that Lemma \ref{three.lagr} does not
always extend to three Lagrangian subspaces.  Suppose that $k =
\mathbb Z/2\mathbb Z$, that $V = k^4$, and that $q = x_1 x_3 + x_2
x_4$.  Let $U$, $U'$ and $U''$ be the Lagrangian subspaces given by
$x_1 = x_2 = 0$, by $x_3 = x_4 = 0$, and by $x_1 + x_3 = x_2 + x_4 =
0$, respectively.  Each subspace is of dimension $2$, and each pair of
subspaces has intersection of dimension $0$.  But there is no
Lagrangian subspace of $V$ which is complementary to $U$, to $U'$ and
to $U''$.

\subsection{Locally alternating maps}
Suppose that $f: \mathcal E \hookrightarrow (\mathcal V,q)$ and $g : 
\mathcal F \hookrightarrow (\mathcal V,q)$ are Lagrangian subbundles of 
a twisted orthogonal bundle.  Consider the composite map 
\begin{equation} 
\label{lambda} 
\lambda:\quad \mathcal E\xrightarrow{\,f\,} \mathcal V 
\xrightarrow{\beta} \mathcal V^*(L) \xrightarrow{g^*} 
\mathcal F^*(L)  
\end{equation} 
where $\beta: \mathcal V \xrightarrow{\cong} \mathcal V^*(L)$ is the
isomorphism induced by the quadratic form $q$.  In the special case of
Lemma \ref{graph.alt}, $\lambda$ is the alternating map $\zeta :
\mathcal F \to \mathcal F^*(L)$.  In the general case, the rank of
$\lambda$ may not be even, and thus $\lambda$ may not be locally
alternating.  But this is the only obstruction: when the rank of
$\lambda$ is even we will show that $\lambda$ is locally alternating,
and we will show how to reduce to the even rank case.

We may assume that $X$ is connected.  Then $\lambda$ is either
everywhere of even rank or everywhere of odd rank because the kernel
of $\lambda(x) : \mathcal E(x) \to \mathcal F^*(L)(x)$ is $\mathcal
E(x) \cap \mathcal F(x)$, which is of constant rank modulo $2$ by
Proposition \ref{constant.parity}.  If $\lambda$ is everywhere of odd
rank, then replace the Lagrangian subbundles $\mathcal E$, $\mathcal
F$ of $\mathcal V$ by the Lagrangian subbundles $\mathcal E_1 :=
\mathcal E \oplus \mathcal O_X$ and $\mathcal F_1 := \mathcal F \oplus
L$ of the orthogonal bundle $\mathcal V_1 := \mathcal V \oplus
\mathcal O_X \oplus L$.  This replaces $\lambda$ by
\begin{equation} 
\label{eq.even.rank} 
\lambda_1: \mathcal E_1 = \mathcal E \oplus \mathcal O_X 
\xrightarrow{\sm{\lambda & 0 \\ 0 & 1}} \mathcal F^*(L) \oplus 
\mathcal O_X = \mathcal F_1^*(L), 
\end{equation} 
The rank of $\lambda_1$ is everywhere even, but its kernel and
cokernel are the same as those of $\lambda$.  Notice also that
$\mathcal E_1(x) \cap \mathcal F_1(x) = \mathcal E(x) \cap \mathcal
F(x)$ for all $x \in X$.  Thus by replacing $\lambda$ by $\lambda_1$
if necessary, we can reduce to the case where the rank is everywhere
even.

\begin{proposition}
\label{prop.loc.alt}
The following are equivalent\textup:

\textup{(a)} The rank of $\lambda$ is even everywhere.

\textup{(b)} $\dim_{k(x)} \bigl[ \mathcal E(x) \cap \mathcal F(x)
\bigr] \equiv \rk(\mathcal E) \pmod 2$ for all $x \in X$.

\textup{(c)} $\lambda$ is locally alternating, i.e.\ there exists a
cover of $X$ by open subsets $U$ and isomorphisms $\iota_U : \mathcal
F \rest U \cong \mathcal E \rest U$ such that the compositions
$\lambda \rest U \circ \iota_U$ are alternating.
\end{proposition}

\begin{proof}
The equivalence of (a) and (b) follows from the fact that of $\ker
\lambda(x) = \mathcal E(x) \cap \mathcal F(x)$.  The implication (c)
$\Rightarrow$ (a) is standard.  To prove (a) $\Rightarrow$ (c), use
Proposition \ref{prop.com.comp} to cover $X$ by open subsets $U$ over
each of which $\mathcal E\rest U, \mathcal F\rest U$ have a common
Lagrangian complement $\mathcal M_U$.  Then by Corollary
\ref{common.compl} there exist alternating maps $\zeta_U: \mathcal F
\rest U \to \mathcal F^*(L) \rest U$ such that $\phi_{\mathcal F\rest
U, \mathcal M_U}(\mathcal E \rest U) = \graph {\zeta_U}$.  Let $\alpha
: \mathcal E \rest U \cong \graph{\zeta_U}$ be the isomorphism induced
by $\phi_{\mathcal F \rest U, \mathcal M_U}$, and let $\pi_1 :
\graph{\zeta_U} \cong \mathcal F \rest U$ and $\pi_2 : \graph{\zeta_U}
\to \mathcal F^*(L) \rest U$ be the two projections from
$\graph{\zeta_U} \subset (\mathcal F \oplus \mathcal F^*(L) ) \rest
U$.
\begin{equation}
\label{made.loc.alt}
\begin{diagram}
\mathcal E \rest U & \rTo ^\cong _\alpha & \graph {\zeta_U} & \rTo
^\cong _{\pi_1} & \mathcal F \rest U \\
& \rdTo < {\lambda \rest U} & \dTo <{\pi_2} & \ldTo >{\zeta_U} \\
&& \mathcal F^*(L) \rest U 
\end{diagram}
\end{equation}
Then $\iota_U := (\pi_1 \circ \alpha)^{-1}$ is an isomorphism such
that $\zeta_U = \lambda \rest U \circ \iota_U$ is alternating.
\end{proof}

\subsection{Independence of the common Lagrangian complement}
Unfortunately, the construction which makes $\lambda$ locally
alternating depends on choices of local common Lagrangian complements.
We now look at what happens if we replace one choice by another.

\begin{lemma}
\label{transform}
Let $\mathcal E, \mathcal F \subset (\mathcal V,q)$ be Lagrangian
subbundles, and let $\mathcal M$ and $\mathcal N$ both be common
Lagrangian complements to $\mathcal E$ and $\mathcal F$.  Suppose that
the map $\phi_{\mathcal F,\mathcal M}$ of \eqref{isom.hyp} sends
$\mathcal E$ to the graph of $\zeta : \mathcal F \to \mathcal F
^*(L)$, and sends $\mathcal N$ to the graph of $h : \mathcal F^*(L)
\to \mathcal F$.  Then

\textup{(a)} The maps $\zeta$ and $h$ are alternating.

\textup{(b)} The map $u := 1 - h \zeta$ and its transpose
$u^* = 1-\zeta h$ are both invertible.

\textup{(c)} The isometry $\phi_{\mathcal F,\mathcal N}$ sends
$\mathcal E$ to the graph of the morphism
\[
\zeta u^{-1} = \left( u^{-1} \right)^* (\zeta - \zeta h
\zeta) u^{-1}.
\]
\end{lemma}

\begin{proof}
Part (a) follows from Corollary \ref{common.compl}.

(b) Since $\mathcal N$ and $\mathcal E$ are complementary, their
images $\graph{\zeta}, \graph{h} \subset \mathcal F \oplus
\mathcal F^*(L)$ are also complementary. Thus the map
\[
\mathcal F \oplus \mathcal F^*(L) \rTo^{\sm{1 & h \\ \zeta & 1}}
\mathcal F \oplus \mathcal F^*(L)
\]  
is an isomorphism.  This is equivalent to the composition
\[
\begin{pmatrix} 1 & -h \\ 0 & 1 \end{pmatrix} 
\begin{pmatrix} 1 & h \\ \zeta & 1 \end{pmatrix} = 
\begin{pmatrix} 1 - h \zeta & 0 \\ \zeta & 1 \end{pmatrix}
\]
being invertible, or to $u = 1 - h \zeta$ being invertible.  It
follows that $u^* = 1 - \zeta h$ is also invertible.

(c) We now have isometries
\begin{equation}
\label{EM.EN}
\mathcal F \oplus \mathcal F^*(L) \lTo^{\phi_{\mathcal F,\mathcal
M}} _\cong \mathcal V \rTo^{\phi_{\mathcal F, \mathcal N}} _\cong
\mathcal F \oplus \mathcal F^*(L).
\end{equation}
The left-to-right composition is the identity on the first summand
$\mathcal F$, and sends $\graph{h}$ (corresponding to $\mathcal N$ on
the left) onto the second summand $\mathcal F^*(L)$ (corresponding to
$\mathcal N$ on the right), compatibly with the hyperbolic quadratic
form on $\mathcal F \oplus \mathcal F^*(L)$.  Consequently, the
left-to-right composition is $\sm{1 & -h \\ 0 & 1}$.  To find the
image of $\mathcal E$ on the right, one first goes to the left (where
its image is $\graph{\zeta}$), and then applies the left-to-right
composition.  Therefore the image of $\mathcal E$ on the right is the
composite image
\begin{diagram}
\mathcal F & \rInto ^{\sm{1 \\ \zeta}} & \mathcal F \oplus \mathcal
F^*(L) & \rTo ^{\sm{1 & -h \\ 0 & 1}} _{\text{left-to-right}} &
\mathcal F \oplus \mathcal F^*(L).\\
\uTo <\cong & \ruInto >{\text{left}\rest {\mathcal E}} \\
\mathcal E
\end{diagram}
The above composite map is $\sm{ 1 & -h \\ 0 & 1 } \sm{ 1 \\ \zeta } =
\sm{ u \\ \zeta } = \sm{ 1 \\ \zeta u^{-1} } u$, so the image of
$\mathcal E$ on the right is $\graph{\zeta u^{-1}}$.  But 
\[
\zeta u^{-1}= \left(u ^{-1} \right)^* \left( u^* \zeta \right)
u^{-1} = \left(u ^{-1} \right)^* \left( \zeta - \zeta h
\zeta \right) u^{-1},
\] 
and this completes the proof. 
\end{proof}

\subsection{Degeneracy Loci}  
Let $\mathcal F$ be a vector bundle of constant rank on a scheme $X$,
let $L$ be a line bundle, and let $\zeta: \mathcal F \to \mathcal F
^*(L)$ be an alternating map.  If $k \geq 0$ is an integer, and if $m
:= \rk(\mathcal F)-2k$, then the degeneracy locus
\begin{equation} 
\label{alt.degen} 
Z_m (\zeta) := \{ x \in X \mid \rk (\zeta(x)) \leq \rk(\mathcal F) - 
m  = 2k \} 
\end{equation} 
has codimension at most $m(m-1)/2$, its ``expected'' value.  The
natural scheme structure on $Z_m (\zeta)$ is defined locally by the
ideal $\Pf_{2k+2}(\zeta)$ generated by the 
$(2k+2) {\times} (2k+2)$
Pfaffians of the alternating map $\zeta$.  These loci have been
studied notably in Harris-Tu \cite{HarrisTu}, and from a different
point of view in Okonek \cite{Okonek} and Walter \cite{Walter}.

We need to extend this notion to the case of a locally alternating
map, as described in Proposition \ref{prop.loc.alt}.  Let $\mathcal E,
\mathcal F \subset (\mathcal V,q)$ be Lagrangian subbundles of a
twisted orthogonal bundle, let $m$ be an integer such that $m \equiv
\dim_{k(x)}\left[ \mathcal E(x) \cap \mathcal F(x) \right] \pmod 2$
for all $x \in X$, and let
\begin{equation} 
\label{lagr.degen} 
Z_{m}(\mathcal E,\mathcal F) := \{ x\in X \mid \dim_{k(x)}[\mathcal 
E(x)\cap \mathcal F(x)]\geq m \}.  
\end{equation} 
The fundamental classes of these loci are discussed in  
Fulton-Pragacz \cite{FP} Chap.\ 6, where they are given as polynomials 
in the Chern classes of $\mathcal E$, $\mathcal F$, and $L$.

We now define a scheme structure on these Lagrangian
degeneracy loci.  Replacing $\lambda$ by the $\lambda_1$ of
\eqref{eq.even.rank} if necessary, we may assume that $\lambda$ is
everywhere of even rank.  By Proposition \ref{prop.loc.alt} $\lambda$
is then locally alternating, i.e.\ there exists a cover of $X$ by open
subsets $U$ and isomorphisms $\iota_U : \mathcal F \rest U \cong
\mathcal E \rest U$ such that the compositions $\zeta_U = \lambda
\rest U \circ \iota_U$ are alternating.
The scheme $Z_{m}(\mathcal E,\mathcal F) \rest U$ is then defined by
\begin{equation} 
\label{scheme.str} 
\mathcal I_{Z_{m} (\mathcal E, \mathcal F)} \rest U := 
\Pf_{\rk(\mathcal E)-m+2}(\zeta_U). 
\end{equation} 
Since $Z_m(\mathcal E, \mathcal F) \rest U$ is the degeneracy locus of
an alternating map, its codimension is at most $m(m -1)/2$.
Now the construction of the maps $\zeta_U$ in Proposition
\ref{prop.loc.alt} depended on the choice of a common Lagrangian
complement to $\mathcal E \rest U$ and $\mathcal F \rest U$.
Nevertheless, the local degeneracy loci $Z_m (\mathcal E, \mathcal F)
\rest U$ are independent of this choice and therefore glue together to
form a scheme $Z_m(\mathcal E, \mathcal F)$ because of Lemma
\ref{transform} and the next lemma.

\begin{lemma}
\label{lambda.beta}
Let $\mathcal F$ be a vector bundle of constant rank and $L$ a line bundle
over a scheme $X$.  If $\zeta : \mathcal F \to \mathcal F^*(L)$ and
$h : \mathcal F^*(L) \to \mathcal F$ are alternating maps such
that $u := 1 - h \zeta$ is invertible, then $\Pf_{2k}
(\zeta) = \Pf_{2k} (\zeta - \zeta h \zeta)$ for
all integers $k$.
\end{lemma}

\begin{proof}
It is enough to prove the lemma in the case where $L = \mathcal O_X$
and where $X$ is universal.  So let $r := \rk (\mathcal F)$, and let
$R := \mathbb Z[X_{ij},Y_{ij}]$ be the polynomial ring in the
independent variables $X_{ij}, Y_{ij}$ ($1 \leq i < j \leq r$).  Let
$\zeta$ and $h$ be the $r {\times} r$ matrices with coefficients in
$R$ given by
\begin{align*}
\zeta_{ij} & := \begin{cases} X_{ij} & \text{if $i < j$,} \\ 0 &
\text{if $i = j$,} \\ -X_{ji} & \text{if $i > j$,} \end{cases} &
h_{ij} & := \begin{cases} Y_{ij} & \text{if $i < j$,} \\ 0 &
\text{if $i = j$,} \\ -Y_{ji} & \text{if $i > j$,} \end{cases}
\end{align*}
let $u := 1 - h \zeta$, let $\delta := \det(u)$, and set
$R_\delta := R[\frac 1 \delta]$. Thus $X := \Spec(R_\delta)$.  
We have to show that the two ideals
\begin{align*}
I & := \Pf_{2k} (\zeta), & J & := \Pf_{2k}(\zeta - \zeta
h \zeta)
\end{align*}
in $R_\delta$ coincide.  However, we may notice the following three
facts:

{\em The ideal $I \subset R_\delta$ is prime.}  This is because the
ideal of $\mathbb Z[X_{ij}]$ generated by the $2k {\times} 2k$
Pfaffians of $\zeta$ is prime (cf.\ Abeasis-Del Fra \cite{Abeasis} \S
3), so its extensions to $R$ (which is a polynomial algebra over over
$\mathbb Z[X_{ij}]$) and to $R_\delta$ are also prime.

{\em There is an involution of $R_\delta$ exchanging $I$ and
$J$.}  Since $R$ is a polynomial algebra over $\mathbb Z$ in variables
which are the entries of $\zeta$ and $h$, one can specify a
morphism $f : R \to R_\delta$ by specifying alternating matrices
$f(\zeta)$ and $f(h)$.  Thus we may define $f$ by
\begin{align*}
f(\zeta) & := \zeta - \zeta h \zeta = \zeta u = u^*
\zeta, & f(h) & := -u^{-1} h \left( u^* \right)^{-1}.
\end{align*}
One computes that $f(u) = u^{-1}$, so $f(\delta)$ is the invertible
element $1/\delta \in R_\delta$.  Hence $f$ extends uniquely to a
morphism $f: R_\delta \to R_\delta$.  One checks that $f(f(\zeta)) =
\zeta$ and that $f(f(h)) = h$, so $f$ is an involution.
Since $f$ exchanges $\zeta$ and $\zeta - \zeta h \zeta$, it
exchanges $I$ and $J$.

{\em The ideals $I$ and $J$ define the same algebraic subset of\/
$\Spec(R_\delta)$.}  This is equivalent to showing that a morphism $g
: R_\delta \to K$ with $K$ a field factors through $R_\delta/I$ if and
only if it factors through $R_\delta / J$.  But giving such a $g$ is
equivalent to giving alternating matrices $g(\zeta)$ and $g(h)$ with
coefficients in $K$ such that $g(u) = 1-g(h) g(\zeta)$ is invertible.
Such a $g$ factors through $R_\delta / I$ if and only if $\rk [
g(\zeta) ] < 2k$, and it factors through $R_\delta / J$ if and only
if
\[
\rk [ g(\zeta - \zeta h \zeta) ] = \rk [ g(\zeta) g(u) ] < 2k.
\]
Since $g(u)$ is invertible, the two conditions are equivalent.

These three facts show that (in the generic case) $I$ and $J$ are
prime ideals defining the same algebraic subsets.  This proves that
$I$ are $J$ are equal in the generic case, and therefore equal in all
cases.  This proves the lemma.
\end{proof}

The definition of $Z_m (\mathcal E, \mathcal F)$ generalizes the
definition of $Z_m(\zeta)$:

\begin{lemma} 
Let $\zeta : \mathcal F \to \mathcal F^*(L)$ be an alternating map
of vector bundles, let $\Gamma_\zeta(\mathcal F) \subset \mathcal F
\oplus\mathcal F^*(L)$ be its graph.  For any $m$ such that $m \equiv
\rk(\mathcal F) \pmod 2$ the degeneracy loci $Z_m(\zeta)$ and
$Z_m(\mathcal F, \Gamma_\zeta(\mathcal F))$
are identical schemes.
\end{lemma} 

Finally, we may verify that our scheme-theoretic degeneracy loci do
not change if we invert the order of our pair of Lagrangian
subbundles, i.e.\ 
\[
Z_m(\mathcal E,\mathcal F) = Z_m(\mathcal F,\mathcal E)
\]
for all $m \equiv \rk(\mathcal E) \pmod 2$.  Essentially, if $\mathcal
M_U$ is a common Lagrangian complement of $\mathcal E \rest U$ and of
$\mathcal F \rest U$ which leads to a $\zeta_U$ as in
\eqref{scheme.str} such that $Z_m(\mathcal E, \mathcal F) \rest U =
Z_m (\zeta_U)$, then a computation similar to Lemma \ref{transform}(c)
leads to a natural identification $Z_m(\mathcal F, \mathcal E) \rest U
= Z_m ( - \zeta_U)$.  We leave the details to the reader.

Our interest in this paper is in the locus $Z := Z_3(\mathcal E,
\mathcal F)$ in the case when it has codimension $3$, the largest
possible (and expected) value.  The fundamental classes computed in
\cite{FultonHirz}, \cite{FultonJDG} and \cite{FP} agree with the
scheme structures introduced here, which coincide with the
scheme structures defined in \cite{DeConciniPragacz}.

All the results of this section have analogues for pairs of
Lagrangian subbundles of twisted symplectic bundles.  Degeneracy loci
for such pairs are the degeneracy loci of locally symmetric maps.  The
symplectic case is slightly simpler than the orthogonal case because
one does not need to worry about the parity of $m$ or of $\dim_{k(x)}
\left[ \mathcal E(x) \cap \mathcal F(x) \right]$ (there is no
symplectic analogue of Proposition \ref{constant.parity}).  The
details are left to the reader.

\section{Lagrangian degeneracy loci are strongly subcanonical} 
\label{Subcanonical}

We now prove the implication (b) $\Rightarrow$ (a) of our main Theorem
\ref{main.result}.
 
\begin{theorem} 
\label{second.v} 
Suppose that $(\mathcal V,q)$ is a twisted orthogonal bundle over a 
locally Noetherian scheme $X$ with values in a line bundle $L$. 
Suppose that $\mathcal E, \mathcal F \subset (\mathcal V,q)$ are 
Lagrangian subbundles such that $\dim_{k(x)} \bigl[ \mathcal E(x) \cap 
\mathcal F(x) \bigr]$ is odd for all $x \in X$.  Write $L_{\mathcal 
E,\mathcal F,\mathcal V} := \det(\mathcal E) \otimes \det(\mathcal F) 
\otimes \det(\mathcal V)^{-1}$. Suppose that the submaximal minors of 
the composite map $\lambda : \mathcal E \hookrightarrow \mathcal V 
\cong \mathcal V^*(L) \twoheadrightarrow \mathcal F^*(L)$ generate an 
ideal sheaf $\mathcal I$ of grade $3$ \textup(the expected 
value\textup). 
 
Then the ideal sheaf of the closed subscheme \textup(cf.\  
\eqref{scheme.str}\textup) 
\[ 
Z ={Z_{3}(\mathcal E,\mathcal F)}=  
\{ x \in X \mid \dim_{k(x)} \left[ \mathcal E(x) \cap \mathcal 
F(x) \right] \geq 3 \}, 
\] 
has grade $3$ and satisfies $\mathcal I_Z^2 = \mathcal I$.  The sheaf 
$\mathcal O_Z$ has locally free resolutions 
\begin{subequations} 
\begin{gather} 
\label{res.1} 
0 \to L_{\mathcal E,\mathcal F,\mathcal V} \to \mathcal E(M) 
\xrightarrow{\,\lambda\,} \mathcal F^*(L\otimes M) \to \mathcal O_X 
\to \mathcal O_Z \to 0, \\ 
\label{res.2} 
0 \to L_{\mathcal E,\mathcal F,\mathcal V} \to \mathcal F(M) 
\xrightarrow{-\lambda^*} \mathcal E^*(L\otimes M) \to \mathcal O_X \to 
\mathcal O_Z \to 0, 
\end{gather} 
\end{subequations} 
with $M$ a line bundle such that $M^{\otimes 2} \cong L_{\mathcal
E,\mathcal F,\mathcal V} \otimes L^{-1}$.  Moreover, the natural
isomorphism between \eqref{res.2} and the dual of \eqref{res.1}
defines an isomorphism
\[ 
\eta: \mathcal O_Z \xrightarrow{\cong} \shExt_{\mathcal O_X}^3 
(\mathcal O_Z, L_{\mathcal E,\mathcal F,\mathcal V}) =: \omega_{Z/X} 
(L_{\mathcal E,\mathcal F,\mathcal V}), 
\] 
with respect to which $Z$ is strongly subcanonical of codimension $3$
in $X$ \textup(cf.\ Definition \ref{strong.def}\textup).
\end{theorem}

\begin{corollary} 
If, in the situation of the theorem, $X$ is locally Gorenstein, then 
so is $Z$, and $\omega_Z \cong \omega_X(L_{\mathcal E,\mathcal 
F,\mathcal V}^{-1}) \rest Z$. 
\end{corollary}

The statement of the theorem remains true even if $X$ is not 
Noetherian, provided one defines grade as in Eagon-Northcott \cite{EN} 
and Northcott \cite{Northcott}.  The only difference in the proofs is 
that one uses the non-Noetherian generalizations of the 
Buchsbaum-Eisenbud structure theorems found in these references.

\begin{proof}[Proof of Theorem \ref{second.v}] 
Let
$f:\ \mathcal E\to \mathcal V$ and 
$g:\ \mathcal F\to \mathcal V$ be the inclusions, and let
 $N = \mathcal E \cap \mathcal F$ be the kernel in the 
natural sequence
\[ 
0 \to N \xrightarrow{\sm{i\\j}} \mathcal E \oplus \mathcal F 
\xrightarrow{\sm{f & -g}} \mathcal V. 
\] 
If $\beta: \mathcal V \xrightarrow{\sim} \mathcal V^*(L)$ is the 
isomorphism induced by the quadratic form $q$, and $\lambda := g^* 
\beta f$, then we get a commutative diagram 
\begin{equation}
\label{pair} 
\begin{diagram}[h=1.5em] 
&& \mathcal E &&\rTo ^\lambda && \mathcal F^*(L) \\ 
& \ruInto^i && \rdInto^f && \ruOnto^{g^*\beta} && \rdTo^{j^*} \\ 
N && && \mathcal V &&&& N^{-1} \otimes L. \\ 
& \rdInto_j && \ruInto^g && \rdOnto^{f^*\beta} && \ruTo_{i^*} \\ 
&& \mathcal F && \rTo_{\lambda^*} && \mathcal E^*(L) 
\end{diagram} 
\end{equation}
Since the diagonals are short exact sequences, the kernels of 
$\lambda$ and of $\lambda^*$ are both equal to $N$.  In addition, 
$f i = g j$.

We claim that $N$ is a line bundle, and that the complexes 
\begin{subequations} 
\begin{gather} 
\label{N} 
0 \to N \xrightarrow{\,i\,} \mathcal E \xrightarrow{\,\lambda\,} 
\mathcal F^*(L) \xrightarrow{j^*} N^{-1}\otimes L \\ 
\label{N2} 
0 \to N \xrightarrow{\,j\,} \mathcal F \xrightarrow{-\lambda^*} 
\mathcal E^*(L) \xrightarrow{i^*} N^{-1}\otimes L  
\end{gather} 
\end{subequations} 
are exact and are locally free resolutions of $\mathcal 
O_Z(N^{-1}\otimes L)$ for the subscheme $Z =Z_{3}(\mathcal E,\mathcal 
F)\subset X$ of grade $3$, with $\mathcal I_Z^2 = \mathcal I$.  We 
will prove these claims locally by making $\lambda$ locally 
alternating and applying the Buchsbaum-Eisenbud structure theorem 
\cite{BE}. 
 
Now the vector bundles $\mathcal E$ and $\mathcal F$ may be of even or
odd rank.  If the rank is even, we use the same trick as in
\eqref{eq.even.rank} and replace $\lambda$ by
\[ 
\mathcal E \oplus \mathcal O_X \xrightarrow{\sm{\lambda & 0 \\ 0 & 1}} 
\mathcal F^*(L) \oplus \mathcal O_X 
\] 
without changing the kernel and cokernel of $\lambda$.  Thus we may 
assume that $\mathcal E$ and $\mathcal F$ are of odd rank. 
 
By hypothesis $\dim_{k(x)} \bigl[ \mathcal E(x) \cap \mathcal F(x)
\bigr]$ is also odd for all $x \in X$.  Therefore, by Proposition
\ref{prop.loc.alt}, $\lambda$ is locally alternating, i.e.\ $X$ is
covered by open subsets $U$ over which there are isomorphisms $\iota_U
: \mathcal F \rest U \cong \mathcal E \rest U$ such that $\zeta_U =
\lambda \rest U \circ \iota_U$ is alternating.  Thus we see that our
complexes \eqref{N} and \eqref{N2} are locally isomorphic to complexes
\begin{equation} 
\label{BE.cpx} 
0 \to N \rest U \xrightarrow{j \rest U} \mathcal F
\xrightarrow{\zeta_U} \mathcal F^*(L) \xrightarrow{j^*\rest U} (N^{-1} 
\otimes L) \rest U 
\end{equation} 
such that $\zeta_U$ is alternating with kernel $j \rest U =
\iota_U^{-1} \circ i \rest U$ in the notation of the diagram
\eqref{pair}.  Now $\mathcal F$ is of odd rank, $\zeta_U$ is
alternating, and the ideal $\mathcal I$ generated by its submaximal
minors is of grade $3$.  So the Buchsbaum-Eisenbud structure theorem
\cite{BE} applies.  Therefore the kernel $N \rest U$ is a line bundle,
the map $j \rest U$ is given by the submaximal Pfaffians of $\zeta_U$,
and the complex \eqref{BE.cpx} is exact and is a resolution of
$\mathcal O_Z(N^{-1}\otimes L) \rest U$.  We can also identify the
ideal sheaf $\mathcal I$ generated by the submaximal minors of
$\lambda$ with $\mathcal I_Z^2$.  This works because on $U$ the sheaf
$\mathcal I_Z \rest U$ is generated by the submaximal Pfaffians $p_1,
\dots , p_n$ of the alternating map $\zeta_U$, while $\mathcal I \rest
U$ is generated by the submaximal minors.  Since the $(i,j)$-th
submaximal minor is $\pm p_i p_j$ (\cite{BE} Appendix), we do indeed
get $\mathcal I_Z^2 \rest U = \mathcal I \rest U$.  This verifies our
claims.
 
Now we set $M := N \otimes L^{-1}$ and twist.  We get two dual 
resolutions 
\begin{subequations} 
\begin{align} 
\label{normal.res} 
0 \to M^{\otimes 2} \otimes L \to \mathcal E(M) & 
\xrightarrow{\,\lambda\,} \mathcal F^*(L\otimes M) \to \mathcal O_X 
\to \mathcal O_Z \to 0, \\  
\label{normal.res2} 
0 \to M^{\otimes 2} \otimes L \to \mathcal F(M) & 
\xrightarrow{-\lambda^*} \mathcal E^*(L\otimes M) \to \mathcal O_X \to 
\mathcal O_Z \to 0. 
\end{align} 
\end{subequations} 
The alternating product of the determinant line bundles in each 
resolution is trivial, and therefore $M^{\otimes 2} \otimes L \cong 
L_{\mathcal E,\mathcal F,\mathcal V}$. 
 
Let us now verify that $Z \subset X$ satisfies conditions
\cond1--\cond4 of Definition \ref{strong.def}.  The duality between
the two resolutions of $\mathcal O_Z$ shows that $Z \subset X$ is
relatively Cohen-Macaulay of codimension $3$.  The duality also
induces an isomorphism $\eta : \mathcal O_Z \xrightarrow{\sim}
\shExt^3_{\mathcal O_X} (\mathcal O_Z, L_{\mathcal E,\mathcal
F,\mathcal V})$, making $Z \subset X$ subcanonical.  Clearly $\mathcal
O_Z$ is of finite local projective dimension.  Moreover, $\eta$ is the
Yoneda extension class of \eqref{normal.res} and is thus the image of
the class in $\Ext^2_{\mathcal O_X}(\mathcal I_Z,M^{\otimes 2} \otimes
L)$ of
\[ 
0 \to M^{\otimes 2} \otimes L \to \mathcal E(M) \to \mathcal 
F^*(L\otimes M) \to \mathcal I_Z \to 0. 
\] 
So $Z \subset X$ is strongly subcanonical with respect to $\eta$. 
\end{proof}

\section{Degeneracy loci for split Lagrangian subbundles}
\label{lagrange} 
 
In this section we discuss Lagrangian degeneracy loci in the ``split''
case where the bundles are $\mathcal E, \mathcal F \subset (\mathcal F
\oplus \mathcal F^*(L),\hyp)$.  We prove the implications (c)
$\Leftrightarrow$ (d) $ \Rightarrow$ (a) of our main Theorem
\ref{main.result}.  We discuss the relation between this case and
Pfaffian subschemes, and we show how the general case can be
transformed into the split case.  In \cite{EPW2} we use this split case
to construct non-Pfaffian subcanonical subschemes of codimension $3$
in $\mathbb P^n$.

\begin{theorem} 
\label{first.v}
Let $\mathcal F$ be a vector bundle of rank $n$ and $L$ a line bundle 
on a locally Noetherian scheme $X$.  Let $\mathcal F \oplus \mathcal 
F^*(L)$ be the hyperbolic twisted orthogonal bundle.  

\textup{(a)}
Suppose that 
\begin{equation} 
\label{inclusion} 
\mathcal E \rInto^{\sm{\psi \\ \phi}} \mathcal F \oplus \mathcal 
F^*(L) 
\end{equation} 
is a Lagrangian subbundle such that $\dim_{k(x)} \left[ \mathcal E(x) 
\cap \mathcal F^*(L)(x) \right]$ is odd for all $x \in X$.  Let 
$L_{\mathcal E,\mathcal F} := \det(\mathcal E) \otimes \det(\mathcal 
F)^{-1}$. 
 
If the sheaf of ideals $\mathcal I$ generated by the submaximal minors 
of $\psi$ is of grade $3$ \textup(the expected value\textup), then the 
ideal sheaf of the closed subscheme  
\[ 
Z ={Z_{3}(\mathcal E,\mathcal F^*(L))}=  
\{ x \in X \mid \dim_{k(x)} \left[ \mathcal E(x) \cap \mathcal 
F^*(L)(x) \right] \geq 3 \}, 
\] 
has grade $3$ and satisfies $\mathcal I_Z^2 = \mathcal I$.  There is a 
commutative diagram with exact rows 
\begin{small} 
\begin{equation}
\label{dual.diag} 
\begin{diagram}
0 & \rTo & L_{\mathcal E,\mathcal F} & \rTo & \mathcal E(M) & \rTo 
^\psi & \mathcal F(M)& \rTo & \mathcal O_X & \rOnto & \mathcal O_Z \\ 
&& \dSame && \dTo <\phi && \dTo >{\phi^*} && \dSame && \dTo >\eta 
<\cong \\ 
0 & \rTo & L_{\mathcal E,\mathcal F} & \rTo & \mathcal F^*(L\otimes M) 
& \rTo ^{-\psi^*} & \mathcal E^*(L\otimes M)& \rTo & \mathcal O_X & 
\rOnto & \shExt^3_{\mathcal O_X} (\mathcal O_Z, L_{\mathcal E,\mathcal 
F}) 
\end{diagram}%
\end{equation}%
\end{small}%
with $M$ a line bundle on $X$ such that $M^{\otimes 2} \cong L^{-1} 
\otimes L_{\mathcal E,\mathcal F}$ and with $\phi^* \psi$ alternating. 
Moreover, $\omega_{Z/X} \cong L^{-1}_{\mathcal E,\mathcal F} \rest Z$, 
and $Z$ is strongly subcanonical of codimension $3$ in $X$ with 
respect to $\eta$. 

\textup{(b)} Conversely, given a subscheme $Z$ with locally free
resolutions as in \eqref{dual.diag} and with $\phi^* \psi$
alternating, then $\mathcal E \subset \mathcal F \oplus \mathcal
F^*(L)$ is a Lagrangian subbundle, and thus $Z \subset X$ is strongly
subcanonical of codimension $3$.
\end{theorem}

\begin{corollary} 
If, in the situation of the theorem, $X$ is locally Gorenstein, then 
so is $Z$, and $\omega_Z \cong \omega_X(L_{\mathcal E,\mathcal 
F}^{-1}) \rest Z$. 
\end{corollary} 
 
\begin{proof}[Proof of Theorem \ref{first.v}] 
  The only things we need to show for part (a) which do not already follow from 
  Theorem \ref{second.v} are that \eqref{dual.diag} commutes and that 
  $\phi^*\psi$ is alternating.  But $\mathcal E$ is a totally isotropic 
  subbundle, so any local section $e \in \Gamma(U,\mathcal E)$ 
  satisfies 
\begin{equation} 
\label{isotropic} 
0 = \hyp(\psi(e)\oplus\phi(e)) = \langle \psi(e),\phi(e) \rangle = 
\langle \phi^*\psi(e),e \rangle. 
\end{equation} 
Thus $\phi^*\psi$ is alternating, and the central part of diagram 
\eqref{dual.diag} commutes.  The rest is easy and left to the reader.
For Part (b) the fact that \eqref{dual.diag} is a quasi-isomorphism implies that $\mathcal E$
is a subbundle of $\mathcal F \oplus \mathcal F^*(L)$.  It is totally isotropic
by the same calculation \eqref{isotropic}.
\end{proof} 
 
\subsection{Pfaffian subschemes} 
Okonek's Pfaffian subschemes \cite{Okonek} are the special case of the 
construction of Theorem \ref{first.v} with $\mathcal E = \mathcal 
F^*(L)$ and $\phi = 1$.  For if 
\[ 
\mathcal E \rInto^{\sm{\psi \\ 1}} \mathcal E^*(L) \oplus \mathcal E 
\] 
is a Lagrangian subbundle, then $\psi$ is alternating by Lemma 
\ref{graph.alt} or by \eqref{isotropic}.  So in this case the two 
resolutions of \eqref{dual.diag} reduce to 
\[ 
0 \to L \otimes M^{\otimes 2} \to \mathcal E(M) \xrightarrow{\psi} 
\mathcal E^*(L \otimes M) \to \mathcal O_X \to \mathcal O_Z \to 0, 
\] 
with $\psi$ alternating.  Thus $Z \subset X$ is one of Okonek's 
Pfaffian subschemes.

\subsection{From non-split to split bundles} 
When the orthogonal bundle $(\mathcal V,q)$ of Theorem \ref{second.v} 
does not split as $\mathcal F \oplus \mathcal F^*(L)$, we cannot fill 
in the diagram \eqref{dual.diag} with direct arrows 
\begin{equation}
\label{dash} 
\begin{diagram}
0 & \rTo & L_{\mathcal E,\mathcal F,\mathcal V} & \rTo ^i & \mathcal 
E(M) & \rTo ^\lambda & \mathcal F^*(L\otimes M) & \rTo ^{j^*} & 
\mathcal O_X & \rOnto & \mathcal O_Z \\  
&& \dSame && \dDotsto && \dDotsto && \dSame && \dTo >\eta 
<\cong \\ 
0 & \rTo & L_{\mathcal E,\mathcal F,\mathcal V} & \rTo ^j & 
\mathcal F(M) & \rTo ^{-\lambda^*} & \mathcal E^*(L\otimes M) & 
\rTo ^{i^*} & \mathcal O_X & \rOnto & \shExt^3_{\mathcal O_X} 
(\mathcal O_Z, L_{\mathcal E,\mathcal G}). 
\end{diagram}
\end{equation} 
Nevertheless, by modifying the orthogonal bundle and its Lagrangian 
subbundles, we can usually realize the same degeneracy locus as a 
Lagrangian degeneracy locus of a split orthogonal bundle.   
 
We know two strategies to accomplish this under different 
hypotheses.  One is to apply the converse structure theorem 
(Theorem \ref{converse}).  The other strategy works whenever 
the quadratic form $q \in \Gamma(X,(S^2 \mathcal V^*)(L))$ is 
the image of an $\alpha \in \Gamma(X,(\mathcal V^* \otimes 
\mathcal V^*)(L))$, for instance if $2 \in \Gamma(X,\mathcal 
O_X)^\times$.  Then the orthogonal direct sum $(\mathcal V,q) 
\perp (\mathcal V,-q)$ is hyperbolic because of the inverse 
isometries 
\begin{diagram} 
(\mathcal V,q) \perp (\mathcal V,-q) & \pile{\rTo ^{\sm{1&1\\ 
\alpha & -\alpha^*}} \\ \lTo _{\sm{\beta^{-1}\alpha^* & 
\beta^{-1} \\ \beta^{-1}\alpha & -\beta^{-1} }}} & \mathcal V 
\oplus \mathcal V^*(L) 
\end{diagram} 
with $\beta := \alpha+\alpha^*$ the nonsingular symmetric bilinear 
form associated to $q$. 
 
The composite map $\mathcal E \oplus \mathcal F \hookrightarrow 
(\mathcal V,q)\perp(\mathcal V,-q) \cong \mathcal V \oplus \mathcal 
V^*(L)$, or more explicitly 
\[ 
\mathcal E \oplus \mathcal F \xrightarrow{\sm{f&g\\ \alpha f & 
-\alpha^* g}} \mathcal V \oplus \mathcal V^*(L) 
\] 
embeds $\mathcal E \oplus \mathcal F$ as a Lagrangian subbundle of the 
hyperbolic bundle $\mathcal V \oplus \mathcal V^*(L)$.  We may then 
fill in the diagram \eqref{dash} with a sequence of quasi-isomorphisms 
going in both directions: 
\begin{diagram} 
0 & \rTo & L_{\mathcal E,\mathcal F,\mathcal V} & \rTo & \mathcal E(M) 
& \rTo ^{g^*\beta f} & \mathcal F^*(L\otimes M) & \rTo & 
\mathcal O_X \\ 
&& \dSame && \uOnto <{\sm{1&0}} && \uOnto >{g^*\beta} && \dSame 
\\ 
0 & \rTo & L_{\mathcal E,\mathcal F,\mathcal V} & \rTo & (\mathcal E 
\oplus \mathcal F)(M) & \rTo ^{\sm{f&g}} & \mathcal V(M) & \rTo & 
\mathcal O_X \\ 
&& \dSame && \dTo <{\sm{\alpha f & -\alpha^* g}} && \dTo 
>{\sm{f^*\alpha^* \\ -g^* \alpha^*}} && \dSame \\ 
0 & \rTo & L_{\mathcal E,\mathcal F,\mathcal V} & \rTo & \mathcal 
V^*(L\otimes M) & \rTo _{\sm{-f^*\\-g^*}} & (\mathcal E^*\oplus 
\mathcal F^*)(L\otimes M) & \rTo & \mathcal O_X \\ 
&& \dSame && \uInto <{\beta g} && \uInto >{\sm{1\\0}} && \dSame 
\\  
0 & \rTo & L_{\mathcal E,\mathcal F,\mathcal V} & \rTo & 
\mathcal F(M) & \rTo _{-f^*\beta g} & \mathcal E^*(L\otimes M) 
& \rTo & \mathcal O_X  
\end{diagram} 

Another way of looking at this is to say: let $\mathcal P_\cpx$ and 
$\mathcal Q_\cpx$ denote the first two lines of the last diagram.  If 
one has $\mathcal V \cong \mathcal F \oplus \mathcal F^*(L)$ as in 
Theorem \ref{first.v}, then the chain map of that theorem is induced 
by a twisted shifted nonsingular quadratic form on the chain complex 
$\mathcal P_\cpx$ given by a chain map $D_2(\mathcal P_\cpx) \to 
L_{\mathcal E,\mathcal F,\mathcal V}[3]$.  In general, no such chain 
map exists, but if we can lift $q$ to $\alpha$ as above, then there is a 
pair of chain maps 
\begin{equation} 
\label{frac.quad} 
D_2(\mathcal P_\cpx) \xleftarrow{\,\sim} D_2(\mathcal Q_\cpx) \to 
L_{\mathcal E,\mathcal F,\mathcal V}[3] 
\end{equation} 
with the first arrow a quasi-isomorphism.  This means that in Theorem 
\ref{second.v}, we are also dealing with a sort of twisted shifted 
nonsingular quadratic form on $\mathcal P_\cpx$, but only in the 
derived category.

\section{Local equations for the degeneracy locus} 
\label{local.eq} 
 
Let $Z \subset X$ be a subcanonical subscheme of codimension $3$ which
is a split Lagrangian degeneracy locus as in diagram
\eqref{dual.diag}.  We give two strategies for computing equations
which define this degeneracy locus locally.  The first is based on the
idea of using a common Lagrangian complement to make $\lambda$
alternating as in Proposition \ref{prop.loc.alt}.  The second is based
on finding standard local forms for a pair of Lagrangian submodules.

\subsection{Strategy 1: Alternating homotopies} 
We start with a pair of Lagrangian subbundles $\mathcal E$ and
$\mathcal F^*(L)$ of a twisted orthogonal bundle.  We may reduce to
the case where the rank of $\mathcal E$ and $\mathcal F$ are odd using
\eqref{eq.even.rank}.  Also since we are working locally, we may
assume that the orthogonal bundle splits as $\mathcal F \oplus
\mathcal F^*(L)$ because over an affine scheme any Lagrangian
subbundle has a Lagrangian complement (e.g.\ Knus \cite{Knus} Remark
I.5.5.4).  Then locally $\mathcal E$ and $\mathcal F^*(L)$ have a
common Lagrangian complement $\mathcal M$ according to Proposition
\ref{prop.com.comp}.  This $\mathcal M$ is necessarily the graph of an
alternating map $h : \mathcal F \to \mathcal F^*(L)$ by Lemma
\ref{graph.alt}.  We use this $h$ as an alternating local homotopy to
transform the commutative diagram on the left below into the one on
the right
\begin{equation} 
\begin{diagram}[inline] 
\label{center} 
\mathcal E & \rTo ^\psi & \mathcal F \\ 
\dTo <\phi & \ldDotsto <h & \dTo >{\phi^*} \\ 
\mathcal F^*(L) & \rTo _{-\psi^*} & \mathcal E^*(L) 
\end{diagram} 
\qquad \quad 
\begin{diagram}[inline] 
\mathcal E & \rTo ^{\psi} & \mathcal F  \\ 
\dTo <{\phi - h \psi} && \dTo >{\phi^* + \psi^*h}  \\ 
\mathcal F^*(L)  & \rTo _{-\psi^*} & \mathcal E^*(L) 
\end{diagram} 
\end{equation} 
Then $\phi -h \psi$ is an isomorphism because it is the projection of 
$\mathcal E$ onto $\mathcal F^*(L)$ along their common complement 
$\mathcal M$.  The dual map $\phi^* + \psi^*h$ is also an isomorphism. 
Thus diagram \eqref{dual.diag}, with a symmetric quasi-isomorphism 
between the resolutions of $\mathcal O_Z$, can be modified locally by 
an alternating homotopy to get a diagram (valid locally) with a 
symmetric isomorphism from the resolution into its dual 
\begin{small} 
\begin{diagram} 
0 & \rTo & L_{\mathcal E,\mathcal F} & \rTo & \mathcal E(M) & \rTo 
^\psi & \mathcal F(M)& \rTo & \mathcal O_X & \rOnto & \mathcal O_Z \\ 
&& \dSame && \dTo <{\phi-h\psi} && \dTo >{\phi^*+\psi^*h} && \dSame && 
\dTo >\eta 
<\cong \\ 
0 & \rTo & L_{\mathcal E,\mathcal F} & \rTo & \mathcal F^*(L\otimes M) 
& \rTo ^{-\psi^*} & \mathcal E^*(L\otimes M)& \rTo & \mathcal O_X & 
\rOnto & \shExt^3_{\mathcal O_X} (\mathcal O_Z, L_{\mathcal E,\mathcal 
F}) 
\end{diagram}%
\end{small}%
Thus locally $\mathcal O_Z$ has a symmetric resolution 
\[ 
0 \to L_{\mathcal E,\mathcal F} \to \mathcal E(M) \xrightarrow{\mu} 
\mathcal E^*(L\otimes M) \to \mathcal O_X \to \mathcal O_Z \to 0 
\] 
where $\mu = \phi^*\psi+\psi^*h\psi$ is alternating (and is 
essentially the map $\mu_{\mathcal M}$ of \eqref{BE.cpx} and 
Proposition \ref{prop.com.comp}).  The submaximal Pfaffians of 
$\mu$ give local equations for the degeneracy locus $Z$. 
 
The choice of another common Lagrangian complement $\mathcal M$ 
gives a different alternating homotopy $h$, and vice-versa.   These
calculations are similar to Lemma \ref{transform} above.

\subsection{Strategy 2: Standard local forms} 
Suppose that $R$ is a commutative local ring with maximal ideal 
$\mathfrak m$ and residue field $k := R/\mathfrak m$.  Let $F$ be a 
free $R$-module of finite rank, and equip $F \oplus F^*$ with the 
hyperbolic quadratic form.  Suppose that $E \subset F \oplus F^*$ is a 
Lagrangian submodule, i.e.\ a totally isotropic direct summand of rank 
equal to that of $F$.  Let $\psi : E \to F$ and $\phi: E \to F^*$ be 
the two components of the inclusion. 
 
\begin{lemma} 
\label{form.for.matrix} 
In the above situation, there exist bases of $E$ and $F$ and a dual 
basis of $F^*$ in which the matrices of $\psi$ and $\phi$ are of the 
form 
\begin{align*} 
\psi & = \begin{pmatrix} \beta & 0 \\ 0 & I \end{pmatrix} & 
\phi & = \begin{pmatrix} I & 0 \\ 0 & \gamma \end{pmatrix} 
\end{align*} 
with the blocks in the two matrices of the same size, and with $\beta$ and 
$\gamma$ alternating. 
\end{lemma} 
 
\begin{proof} 
We begin by choosing bases for $E$ and $F$ and the dual basis of 
$F^*$, so that we can treat $\psi$ and $\phi$ as matrices.  Since 
$E$ is a direct summand of $F \oplus F^*$, the columns of the total 
matrix $\sm{\psi \\ \phi}$ are linearly independent even modulo 
$\mathfrak m$.  Moreover, by \eqref{isotropic} above $\phi^*\psi$ is 
an alternating matrix because $E \subset F \oplus F^*$ is a totally 
isotropic submodule. 
 
We now begin a series of row and column operations on $\psi$ and 
$\phi$ which will put them into the required form.  The column 
operations (resp.\ row operations) correspond to changes of basis of 
$E$ (resp.\ of $F$ and $F^*$) and to the action of invertible matrices 
$P$ (resp.\ $Q$) on $\psi$ and $\phi$ via $\psi \rightsquigarrow 
Q^{-1}\psi P$ and $\phi \rightsquigarrow Q^* \phi P$. 
 
Choose a maximal invertible minor of $\phi$.  After row and column 
operations, we can assume that the corresponding submatrix is an 
identity block lying in the upper left corner of $\phi$ and that the 
blocks below and to the right of it are $0$.  Thus we can assume that 
\begin{align*} 
\phi & = \begin{pmatrix} I & 0 \\ 0 & \delta \end{pmatrix} & 
\psi & = \begin{pmatrix} \psi_{11} & \psi_{12} \\ \psi_{21} & 
\psi_{22} \end{pmatrix} 
\end{align*} 
where the blocks of the two matrices are of the same size, the 
on-diagonal blocks are square, and the coefficients of $\delta$ lie in 
$\mathfrak m$.  Since  
\[ 
\phi^* \psi = \begin{pmatrix} \psi_{11} & \psi_{12} \\ \delta^* 
\psi_{21} & \delta^* \psi_{22} \end{pmatrix} 
\] 
is alternating, we see that all the coefficients of $\psi_{12}$ also 
lie in $\mathfrak m$.  Hence all the coefficients in the last block of 
columns of $\sm{\psi \\ \phi}$ lie in $\mathfrak m$ except those in 
$\psi_{22}$.  Since these columns must be linearly independent modulo 
$\mathfrak m$, it follows that $\psi_{22}$ must be invertible. 
Applying a new set of column operations to $\phi$ and $\psi$, we may 
assume that $\psi = \sm{\psi_{11} & \epsilon \\ \psi_{21} & I}$ and 
that $\phi = \sm{I & 0 \\ 0 & \gamma}$.  Moreover, $\phi^* \psi$ 
remains alternating, which actually means that $\psi_{11}$ and 
$\gamma$ are alternating, and $\epsilon = -\psi_{21}^* \gamma$.  A 
final set of row and column operations using the matrices $Q = \sm{I & 
\psi_{21}^* \gamma \\ 0 & I}$ and $P = \sm{I & 0 \\ -\psi_{21} & I}$ 
puts $\psi$ and $\phi$ into the form required by the lemma. 
\end{proof}

\begin{corollary} 
Let $R$ be a commutative local ring with residue field $k$, let $F$ be
a free $R$-module of odd rank, and let $E \subset F \oplus F^*$ be a
Lagrangian submodule such that $\dim_k \left[ (E\otimes k) \cap (F
\otimes k) \right]$ is odd.  Let $\psi: E \to F$ and $\phi: E \to F$
be the two components of the inclusion.
 
\textup{(a)} The determinant of $\phi$ is of the form $\det \phi = 
af^2$ with $a$ invertible. 
 
\textup{(b)} If $\det(\phi)$ is not a zero-divisor, and if $\psi$ 
degenerates along an ideal $I$ of height and grade $3$ \textup(as 
expected\textup), then this ideal is  $I = 
(\Pf(\phi^*\psi) : f)$, where $\Pf(\phi^*\psi)$ is the ideal generated by 
the submaximal Pfaffians of $\phi^* \psi$, and where $f$ is as in part 
\textup{(a)}. 
 
\end{corollary} 
 
\begin{proof} 
(a) We put the matrices of $\phi$ and of $\psi$ in the special form of 
Lemma \ref{form.for.matrix}, and we set $f := \Pf(\gamma)$.  The 
determinant of the matrix of $\phi$ is then $f^2$.  Consequently the 
determinant of the matrix of $\phi$ with respect to any bases of $E$ 
and $F$ is of the form $af^2$, with $a$ an invertible element of $R$ 
coming from the determinants of the change-of-basis matrices. 
 
(b) Using the special forms for $\phi$ and $\psi$ given in Lemma 
\ref{form.for.matrix}, we find that $I$ is generated by the submaximal 
Pfaffians $p_1, \dots, p_{2s+1}$ of $\beta$, while the ideal 
$\Pf(\phi^*\psi)$ is generated by $fp_1, \dots , fp_{2s+1}$.  Since we 
suppose $\det(\phi)$ and therefore $f$ are not zero-divisors, this 
gives (b). 
\end{proof}

\section{Subcanonical subschemes are Lagrangian degeneracy loci} 
\label{converse.theorem} 
 
In this section we prove the implication (a) $\Rightarrow$ (d) of our
main Theorem \ref{main.result}.  Taken together with Theorems
\ref{second.v} and \ref{first.v}, this proves the main theorem,
because the implication (c) $\Rightarrow$ (b) is trivial.

\begin{theorem} 
\label{converse} 
Let $A$ be a Noetherian ring, and $X \subset \mathbb P^N_A$ a locally
closed subscheme.  If $Z\subset X$ is a codimension $3$ strongly
subcanonical subscheme \textup(cf.\ Definition
\ref{strong.def}\textup), then there exist vector bundles $\mathcal E$
and $\mathcal G$, a line bundle $L$ on $X$, and an embedding of
$\mathcal E$ as a Lagrangian subbundle of the twisted hyperbolic
bundle $\mathcal G\oplus \mathcal G^*(L)$ such that $Z=Z_3(\mathcal E,
\mathcal G^*(L))$ and ${\mathcal O}_Z$ has symmetrically
quasi-isomorphic locally free resolutions
\begin{equation}
\label{converse.diag} 
\begin{diagram}
0 & \rTo & L & \rTo & \mathcal E & \rTo ^\psi & \mathcal G & \rTo & 
\mathcal O_X & \rTo & \mathcal O_Z & \rTo & 0 \\ 
&& \dSame && \dTo <\phi && \dTo >{\phi^*} && \dSame && \dTo >\eta 
<\cong \\ 
0 & \rTo & L & \rTo & \mathcal G^*(L) & \rTo ^{-\psi^*} & \mathcal 
E^*(L)& \rTo & \mathcal O_X & \rTo & \shExt^3_{\mathcal O_X} (\mathcal 
O_Z, L) & \rTo & 0 
\end{diagram} 
\end{equation}
with $\phi^*\psi : \mathcal E \to \mathcal E^*(L)$ an alternating map.
\end{theorem}

We will need the following two lemmas in the proof of the theorem. 
 
\begin{lemma} 
\label{quasiproj} 
Let $A$ be a Noetherian ring, let $X \subset \mathbb P^N_A$ be a 
locally closed subscheme, let $\mathcal F,\mathcal G$ be coherent 
sheaves on $X$, let $\mathcal M$ be a vector bundle on $X$, and let $p 
> 0$.  
 
\textup{(a)} If $\xi \in \Ext^p_{\mathcal O_X}(\mathcal F,\mathcal 
G)$, then there exists a vector bundle $\mathcal E$ on $X$ and a 
surjection $f: \mathcal E \twoheadrightarrow \mathcal F$ such that the 
pullback class $f^*\xi \in \Ext^p_{\mathcal O_X}(\mathcal E,\mathcal 
G)$ vanishes. 
 
\textup{(b)} If $\zeta \in \Ext^p_{\mathcal O_X}(\Lambda^2 \mathcal M, 
\mathcal G)$, then there exists a surjection of vector bundles 
$\mathcal P \twoheadrightarrow \mathcal M$ such that the pullback of 
$\zeta$ to $\Ext^p_{\mathcal O_X}(\Lambda^2 \mathcal P, \mathcal G)$ 
vanishes. 
\end{lemma} 
 
\begin{proof} 
(a) Extending $\mathcal F$ to a coherent sheaf on the closure 
$\overline X \subset \mathbb P^N_A$ and then applying Serre's Theorem 
A, we see that there exists a surjection of the form $g: \mathcal 
O_X(-n)^{r} \twoheadrightarrow \mathcal F$.  Pulling back gives us a 
class $g^*\xi \in H^p(X,\mathcal G(n))^r$. 
 
Let $R$ be the homogeneous coordinate ring of $\overline X$, and let 
$I \subset R$ be the homogeneous ideal of strictly positive degree 
elements vanishing on the closed subset $\overline X \setminus X$. 
Extend $\mathcal G$ to a coherent sheaf on $\overline X$, and let $G$ 
be a finitely generated graded $R$-module whose associated sheaf is 
this extension of $\mathcal G$.  Then $H^p_*(X, \mathcal G)^r \cong 
H^{p+1}_I(G)^r$.  Consequently, $g^*\xi$, as a member of a local 
cohomology module, is annihilated by some power $I^m$ of $I$.  A 
finite set of homogeneous generators of $I^m$ gives surjections 
$\bigoplus_i R(-a_i) \twoheadrightarrow I^m$ and $\bigoplus_i \mathcal 
O_X(-a_i) \twoheadrightarrow \mathcal O_X$, such that the pullback of 
$g^*\xi$ along the induced map $\bigoplus_i \mathcal O_X(-a_i-n)^r 
\twoheadrightarrow \mathcal O_X(-n)^r$ vanishes. 
 
(b) For the same reasons as in part (a), $\zeta$ is killed by some 
power of $I$.  For convenience we assume that the same $I^m$ as in 
part (a) kills $\zeta$.  Let $\bigoplus_i \mathcal O_X(-a_i) 
\twoheadrightarrow \mathcal O_X$ be the surjection used in part (a). 
Then the surjection $\bigoplus_i \mathcal M(-a_i) \twoheadrightarrow 
\mathcal M$ kills $\zeta$ because the exterior square factors as 
$\Lambda^2 \bigl( \bigoplus_i \mathcal M(-a_i) \bigr) 
\twoheadrightarrow \bigoplus_{i\leq j} (\Lambda^2\mathcal M) 
(-a_i-a_j) \twoheadrightarrow \Lambda^2 \mathcal M$. 
\end{proof}

\begin{lemma}[Serre; \cite{OSS} Lemma 5.1.2] 
\label{Serre} 
Let $A$ be a Noetherian local ring and $M$ a finitely generated 
$A$-module of projective dimension at most $1$.  Suppose that $\zeta 
\in \Ext^1_A(M,A)$ corresponds to the extension 
\( 
0 \to A \to N \to M \to 0. 
\) 
Then $N$ is a free $A$-module if and only if $\zeta$ generates the 
$A$-module $\Ext^1_A(M,A)$. 
\end{lemma}

\begin{proof}[Proof of Theorem \ref{converse}] 
By hypothesis, $\eta$ lifts to a class in $\Ext^2_{\mathcal 
O_X}(\mathcal I_Z,L)$.  By Lemma \ref{quasiproj}(a) there exists a 
vector bundle $\mathcal M$ and a surjection and kernel 
$ 
0 \to \mathcal K \to \mathcal M \to \mathcal I_Z \to 0 
$ 
such that $\eta$ lifts further to a class $\zeta \in \Ext^1_{\mathcal 
O_X}(\mathcal K,L)$.  This defines an extension $0 \to L \to \mathcal 
E \to \mathcal K \to 0$.  Attaching these extensions gives an acyclic 
complex 
\begin{equation} 
\label{long} 
0 \to L \to \mathcal E \to \mathcal M \to \mathcal O_X \to \mathcal 
O_Z \to 0. 
\end{equation} 

We claim that $\mathcal E$ is locally free.  Our reasoning is as 
follows.  Since the local projective dimension of $\mathcal O_Z$ is at 
most $3$, the local projective dimension of $\mathcal K$ is at most 
$1$.  By Lemma \ref{Serre}, $\mathcal E$ will be locally free if 
$\zeta$ generates the sheaf $\shExt^1_{\mathcal O_X} (\mathcal K,L)$. 
Moreover the sheaves $\mathcal O_Z$, $\shExt^3_{\mathcal O_X} 
(\mathcal O_Z,L))$, and $\shExt^1_{\mathcal O_X} (\mathcal K,L)$ are 
all isomorphic, and their respective global sections $1$, $\eta$, and 
$\zeta$ correspond under these isomorphisms.    
\begin{diagram} 
\zeta \in \Ext^1_{\mathcal O_X} (\mathcal K, L) & \rTo & H^0 
(\shExt^1_{\mathcal O_X} (\mathcal K,L)) \\ 
\dTo && \dTo >\cong \\ 
\eta \in \Ext^3_{\mathcal O_X} (\mathcal O_Z, L) & \rTo ^\cong & H^0 
(\shExt^3_{\mathcal O_X} (\mathcal O_Z,L)) & \lTo ^\cong & 
H^0(\mathcal O_Z) \ni 1 
\end{diagram} 
Since $1$ generates $\mathcal O_Z$, the section $\zeta$ generates 
$\shExt^1_{\mathcal O_X} (\mathcal K,L)$.  Thus $\mathcal E$ is 
locally free.

The complex 
\begin{equation} 
\label{A.cpx} 
\mathcal A_\cpx : \qquad 0 \to L \to \mathcal E \to \mathcal M \to 
\mathcal O_X \to 0 
\end{equation} 
is now a locally free resolution of $\mathcal O_Z$.  As in 
Buchsbaum-Eisenbud \cite{BE} and Walter \cite{Walter}, we try to make 
this into a commutative associative differential graded algebra 
resolution of $\mathcal O_Z$ by constructing a map $D_2(\mathcal 
A_\cpx) \to \mathcal A_\cpx$ from the divided square covering the 
identity in degree $0$: 
\begin{small}%
\begin{equation}
\label{dashes} 
\begin{diagram}
\dotsb & \rTo & \mathcal M(L) \oplus D_2\mathcal E & \rTo & L \oplus 
(\mathcal E \otimes \mathcal M) & \rTo & \mathcal E \oplus 
\Lambda^2\mathcal M & \rTo & \mathcal M & \rTo & \mathcal O_X \\ 
&& \dDotsto && \dDotsto && \dDotsto && \dSame && \dSame \\ 
\dotsb & \rTo & 0 & \rTo & L & \rTo & \mathcal E & \rTo & \mathcal M & 
\rTo & \mathcal O_X 
\end{diagram}
\end{equation}%
\end{small}%
Now $\Lambda^2\mathcal M$ maps into the kernel $\mathcal K$ of 
$\mathcal M \to \mathcal O_X$.  Hence the first problem in filling in 
the dotted arrows above is to carry out a lifting 
\begin{diagram}[PS] 
&&&&&& \Lambda^2 \mathcal M \\ 
&&&&& \ldDotsto & \dTo \\ 
0 & \rTo & L & \rTo & \mathcal E & \rTo & \mathcal K & \rTo & 0 
\end{diagram} 
The obstruction to carrying out the lifting is a class $\zeta \in 
\Ext^1_{\mathcal O_X}(\Lambda^2 \mathcal M, L)$.  There is no reason 
for this class to vanish.  So the liftings sought in \eqref{dashes} 
need not exist.  But there is a way around this. 
 
By Lemma \ref{quasiproj}(b) there is a surjection from another vector 
bundle $\mathcal G \twoheadrightarrow \mathcal M$ such that the 
pullback of $\zeta$ to $\Ext^1_{\mathcal O_X}(\Lambda^2 \mathcal G,L)$ 
vanishes.  We now redo the construction of the complex and get 
commutative diagrams with exact rows and columns 
\[ 
\begin{diagram}[size=1.5em] 
&&&& \mathcal R & \rSame & \mathcal R \\ 
&&&& \dInto && \dInto \\ 
0 & \rTo & L & \rTo & \mathcal F & \rTo & \mathcal K' & \rTo & 0 \\ 
&& \dSame && \dOnto & \square & \dOnto \\ 
0 & \rTo & L & \rTo & \mathcal E & \rTo & \mathcal K & \rTo & 0 
\end{diagram} 
\qquad \qquad 
\begin{diagram}[size=1.5em] 
&& \mathcal R & \rSame & \mathcal R \\ 
&& \dInto && \dInto \\ 
0 & \rTo & \mathcal K' & \rTo & \mathcal G & \rTo & \mathcal I_Z & 
\rTo & 0 \\ 
&& \dOnto && \dOnto && \dSame \\ 
0 & \rTo & \mathcal K & \rTo & \mathcal M & \rTo & \mathcal I_Z & \rTo 
& 0  
\end{diagram} 
\] 
This allows us to construct a new complex 
\[ 
\mathcal B_\cpx : \qquad 0 \to L \to \mathcal F \xrightarrow{\psi} 
\mathcal G \to \mathcal O_X \to 0. 
\] 
One sees easily that $\mathcal R$ and therefore $\mathcal F$ are also 
vector bundles.  But this time, the composite map $\Lambda^2 \mathcal 
G \to \mathcal K' \to \mathcal K$, lifts to $\mathcal E$ since the 
obstruction is the class in $\Ext^1_{\mathcal O_X}(\Lambda^2 \mathcal 
G,L)$ which we got to vanish using Lemma \ref{quasiproj}(b).  Since 
the square marked with the $\square$ is cartesian, we get a lifting 
$\Lambda^2 \mathcal G \to \mathcal F$.  The other liftings  
\begin{equation}
\label{B.alg}
\begin{diagram}
\dotsb & \rTo & \mathcal G(L) \oplus D_2\mathcal F & \rTo & L \oplus 
(\mathcal F \otimes \mathcal G) & \rTo & \mathcal F \oplus 
\Lambda^2\mathcal G & \rTo & \mathcal G & \rTo & \mathcal O_X \\ 
&& \dDotsto && \dDotsto && \dDotsto && \dSame && \dSame \\ 
\dotsb & \rTo & 0 & \rTo & L & \rTo & \mathcal F & \rTo & \mathcal G & 
\rTo & \mathcal O_X 
\end{diagram} 
\end{equation}
now occur automatically.  We therefore get a chain map $D_2 \mathcal 
B_\cpx \to \mathcal B_\cpx$ which makes $\mathcal B_\cpx$ into a 
commutative associative differential graded algebra with divided 
powers. 
 
We now claim that having this differential graded algebra structure 
gives us all the properties we want and puts us into the situation of 
Theorem \ref{first.v}.  Indeed, as in Buchsbaum-Eisenbud \cite{BE}, 
the multiplication gives pairings $\mathcal B_i \otimes \mathcal 
B_{3-i} \to \mathcal B_3 = L$, and therefore maps $\mathcal B_i \to 
\mathcal B_{3-i}^*(L)$.  These maps are compatible with the 
differential, and as a result, the following diagram commutes: 
\begin{diagram} 
0 & \rTo & L & \rTo & \mathcal F & \rTo^\psi & \mathcal G & \rTo & 
\mathcal O_X & \rTo & \mathcal O_Z & \rTo & 0 \\ 
&& \dSame && \dTo <\phi && \dTo >{\phi^*} && \dSame && \dTo >\eta 
<\simeq \\ 
0 & \rTo & L & \rTo & \mathcal G^*(L) & \rTo_{-\psi^*} & \mathcal 
F^*(L) & \rTo & \mathcal O_X & \rTo & \shExt^3_{\mathcal O_X}(\mathcal 
O_Z,L) & \rTo & 0 
\end{diagram} 
The top row is exact by construction, and the bottom row is exact 
because it is the dual of the top row which is a resolution of a sheaf 
of grade $3$.  Since $\eta$ is an isomorphism, one sees that 
\[ 
0 \to \mathcal F \xrightarrow{\sm{\psi \\ \phi}} \mathcal G \oplus 
\mathcal G^*(L) \xrightarrow{\sm{\phi^* & \psi^*}} \mathcal F^*(L) \to 
0 
\] 
is exact.  Thus $\mathcal F$ embeds in $\mathcal G \oplus \mathcal 
G^*(L)$ as a subbundle which is totally isotropic for the hyperbolic 
symmetric bilinear form on $\mathcal G \oplus \mathcal G^*(L)$.  The 
subbundle $\mathcal F$ is even totally isotropic for the hyperbolic 
quadratic form, since the restriction of this form to local sections 
of $\mathcal F$ is the function $e \mapsto \langle \phi(e), \psi(e) 
\rangle$, and this function vanishes because the composite map from 
diagram \eqref{B.alg} 
\begin{diagram} 
D_2 \mathcal F & \rTo & \mathcal F \otimes \mathcal G & \rTo & L \\ 
f \otimes f & \rMapsto & f \otimes \psi(f) & \rMapsto & \langle 
\phi(f), \psi(f) \rangle  
\end{diagram} 
factors through $0$ and hence vanishes identically.  Thus $\mathcal F$ 
is a Lagrangian subbundle of $\mathcal G \oplus \mathcal G^*(L)$.
This completes the proof.
\end{proof}

\section{Points in $\mathbb P^3$} 
\label{sect.points} 
 
In this and the following section we discuss several classes of
examples which satisfy some or all of the conditions \cond1-\cond4 of
the definition of a strongly subcanonical subscheme and thus Theorems
\ref{main.result} and \ref{converse} may apply.  Additional geometric
applications and examples can be found in our paper \cite{EPW2}.
 
Okonek \cite{Okonek}, p.\ 429, has shown that any reduced set of
points in $\mathbb P^3$ is Pfaffian.  By carefully analyzing the
constructions of Theorem \ref{converse}, we will describe Pfaffian
resolutions of locally Gorenstein zero-dimensional subschemes in
$\mathbb P^3$ (see Remark \ref{classify}).
 
For a locally Gorenstein zero-dimensional subscheme $Z \subset \mathbb
P^3_k$ over a field $k$, there are many isomorphisms $\eta: \mathcal
O_Z \xrightarrow{\sim} \omega_Z(t)$.  Which triples
$(Z,\omega_{\mathbb P^3}(t),\eta)$ satisfy all the conditions of
Definition \ref{strong.def}, and which do not?  In particular (and
this is the only condition which causes trouble), when does the image
of $\eta$ in $H^3(\mathbb P^3, \omega_{\mathbb P^3}(t))$ vanish?
 
We will use the following notation, see for instance \cite{EP}.   
Let $I \subset R := 
k[x_0,x_1,x_2,x_3]$ be the homogeneous ideal of $Z$, let $A := R/I$ be 
its homogeneous coordinate ring, and let $\omega_A := 
\Ext^3_R(A,R(-4))$ be its canonical module.  Note that $\eta \in 
H^0_*(\omega_Z) \supset \omega_A$.  Also if $M$ is a graded $R$-module, 
then let $M'$ be its dual as a graded $k$-vector space, endowed with 
the natural dual $R$-module structure. 
 
\begin{proposition} 
\label{points} 
Let $Z \subset \mathbb P^3$ be a locally Gorenstein subscheme of
dimension zero, and $\eta : \mathcal O_Z \xrightarrow{\sim}
\omega_Z(t)$ an isomorphism.  Then the triple $(Z,\omega_{\mathbb
P^3}(t),\eta)$ is subcanonical and satisfies conditions \cond1-\cond3
of Definition \ref{strong.def}, and it satisfies condition \cond4 if
and only if $\eta \in \omega_A$.
\end{proposition}

\begin{proof} 
The map $\eta: \mathcal O_Z \to \omega_Z(t)$ may be identified with an
element of $\Ext^3_{\mathcal O_{\mathbb P^3}} (\mathcal O_Z,
\omega_{\mathbb P^3}(t)) \cong H^0(\mathcal O_Z(-t))'$.  The subscheme
$Z\subset \mathbb P^3$ satisfies condition \cond4 for $\eta$ if and
only if $\eta$ is in the image of $\Ext^2_{\mathcal O_{\mathbb P^3}}
(\mathcal I_Z, \omega_{\mathbb P^3}(t)) \cong H^1(\mathcal I_Z(-t))'$.
   
Local duality and Serre duality give identifications 
\[ 
\omega_A := \Ext^3_R(A,R(-4)) \cong H^1_{\mathfrak m}(A)' \cong 
H^1_*(\mathcal I_Z)' 
\] 
and $H^0_*(\omega_Z) \cong H^0_*(\mathcal O_Z)'$ which are compatible
with the inclusions.  So $\eta$ satisfies condition \cond4 of
Definition \ref{strong.def} if and only if $\eta \in \omega_A$.
\end{proof}

\begin{theorem} 
\label{pfaff.points} 
Let $Z \subset \mathbb P^3$ be a locally Gorenstein subscheme of 
dimension $0$, and let $\eta \in H^0(\omega_Z(t))$.  Suppose that 
\textup{(a)} $\eta$ generates the sheaf $\omega_Z$, \textup{(b)} $\eta 
\in \omega_A$, and \textup{(c)} if $t=-2\ell$ is even, then the 
following nondegenerate symmetric bilinear form on $H^0(\mathcal 
O_Z(\ell))$ is metabolic \textup(i.e.\ contains a Lagrangian 
subspace\textup{):} 
\begin{equation} 
\label{bilinear} 
H^0(\mathcal O_Z(\ell)) \times H^0(\mathcal O_Z(\ell)) \to 
H^0(\mathcal O_Z(2\ell)) \xrightarrow{\,\eta\,} H^0(\omega_Z) 
\xrightarrow{\tr} k. 
\end{equation} 
Then there exists a locally free resolution 
\begin{equation} 
\label{pfaff.resol} 
0 \to \mathcal O_{\mathbb P^3}(t-4) \to \mathcal F^*(t-4) 
\xrightarrow{\,\psi\,} \mathcal F \to \mathcal O_{\mathbb P^3} \to 
\mathcal O_Z \to 0 
\end{equation} 
with $\psi$ alternating and $\mathcal I_Z$ generated by the submaximal 
Pfaffians of $\psi$ and such that the Yoneda extension class of 
\eqref{pfaff.resol} is $\eta \in \Ext^3_{\mathcal O_{\mathbb P^3}} 
(\mathcal O_Z, \mathcal O_{\mathbb P^3}(t-4)) \cong H^0(\omega_Z(t))$. 
 
Conversely, if there exists a locally free resolution of $\mathcal 
O_Z$ as in \eqref{pfaff.resol} with $\psi$ alternating, then its 
Yoneda extension class $\eta$ satisfies conditions \textup{(a), (b)} 
and \textup{(c)} above. 
\end{theorem}

In order for the symmetric bilinear form \eqref{bilinear} to be 
metabolic, it is necessary for $\deg(Z)$ to be even.  If the base 
field $k$ is closed under square roots, this is also sufficient. 
 
In any case, the conditions of the theorem always hold if $t$ is large 
and odd and $\eta$ is general.  This proves the following result, 
which was proven for reduced sets of points by Okonek \cite{Okonek}, 
p.\ 429.  
 
\begin{corollary} 
A zero-dimensional subscheme of $\mathbb P^3$ is Pfaffian if and only 
if it is locally Gorenstein. 
\end{corollary}

\begin{proof}[Proof of Theorem \ref{pfaff.points}] 
We show how to start the proof off.  But we will stop when we reach 
the point where it becomes identical to the proof of the main result 
of \cite{Walter}. 
 
Suppose that $Z,t,\eta$ satisfy conditions (a), (b), and (c) of the 
theorem.  Condition (a) implies that the map $\eta: \mathcal O_Z \to 
\omega_Z(t)$ is an isomorphism.  So $\eta$ and Serre duality induce a 
symmetric perfect pairing  
\begin{equation} 
\label{sym.pairing} 
H^0_*(\mathcal O_Z) \times H^0_*(\mathcal O_Z) \xrightarrow{\mult} 
H^0_*(\mathcal O_Z) \xrightarrow{\,\eta\,} H^0_*(\omega_Z(t)) 
\xrightarrow{\tr} k(t) 
\end{equation} 
which pairs $H^0(\mathcal O_Z(n))$ with $H^0(\mathcal O_Z(-n-t))$ for 
all $n$. 
 
Condition (c) implies that $H^0_*(\mathcal O_Z)$ contains a Lagrangian 
submodule $M$ for this symmetric perfect pairing.  Indeed if $t$ is 
odd, one can pick $M := \bigoplus_{n > -t/2} H^0(\mathcal O_Z(n))$. 
If $t$ is even, then there exists a Lagrangian subspace $W \subset 
H^0(\mathcal O_Z(-t/2))$, and one can pick $M := W \oplus \bigoplus_{n 
> -t/2} H^0(\mathcal O_Z(n))$. 
 
The two submodules $A \subset H^0_*(\mathcal 
O_Z)$ and $\omega_A \subset H^0_*(\omega_Z)$ are orthogonal 
complements of each other under the Serre duality pairing; see for 
example \cite{EP}.  Hence 
condition (b), that $\eta \in \omega_A$, implies that $\eta A \subset 
\omega_A$ and therefore that $A = \omega_A^\perp \subset (\eta 
A)^\perp$.  Now the orthogonal complement of $\eta A \subset 
H^0_*(\omega_Z)$ under the Serre duality pairing corresponds to the 
orthogonal complement of $A \subset H^0_*(\mathcal O_Z)$ under our 
pairing \eqref{sym.pairing}.  So condition (b) implies that $A \subset 
A^\perp$.  In other words $A \subset H^0_*(\mathcal O_Z)$ is 
sub-Lagrangian. 
 
It now follows that there exists a Lagrangian submodule $L$ such that 
$0 \subset A \subset L = L^\perp \subset A^\perp \subset 
H^0_*(\mathcal O_Z)$.  For instance, pick $L := A + (M \cap A^\perp)$ 
(cf.\ Knus \cite{Knus} Lemma I.6.1.2). 
 
One easily checks that $A_n = (A^\perp)_n = H^0(\mathcal O_Z(n))$ for 
$n \gg 0$, and that $A_n = (A^\perp)_n = 0$ for $n \ll 0$. 
Consequently $A^\perp/A$ is of finite length.  It has an induced 
nondegenerate symmetric bilinear form, and it has a Lagrangian 
submodule $L/A$. 
 
We now claim that we can construct a locally free resolution  
\[ 
0 \to \mathcal O_{\mathbb P^3}(t-4) \xrightarrow{\,\alpha\,} \mathcal 
F^*(t-4) \xrightarrow{\,\psi\,} \mathcal F \xrightarrow{\,\beta\,} 
\mathcal O_{\mathbb P^3} \to \mathcal O_Z \to 0 
\] 
with $\psi$ alternating and such that $H^1_*(\mathcal F) \cong L/A$, 
and $H^2_*(\mathcal F) = 0$.  Moreover, $\beta$ induces a surjection 
$H^0_*(\mathcal F) \twoheadrightarrow H^0_*(\mathcal I_Z)$.  Different 
pieces of the resolution contribute different pieces of the cohomology 
module $H^0_*(\mathcal O_Z)$.  The submodule $A$ is contributed by 
$\coker H^0_*(\beta)$; the piece $L/A$ by $H^1_*(\mathcal F)$; the 
piece $A^\perp/L$ by $H^2_*(\mathcal F^*(t-4))$; and the piece 
$H^0_*(\mathcal O_Z)/A^\perp$ is contributed by $\ker H^3_*(\alpha)$. 
 
The construction of this resolution and the verification of its 
properties can be done using the Horrocks correspondence by the same 
method as in Walter \cite{Walter}.  It is quite long and we omit the 
details. 
\end{proof} 
 
\begin{remark} 
\label{classify} 
The graded module $A^\perp / A$ above can be thought of as the 
``intermediate cohomology'' or {\em deficiency module} of 
$(Z,t,\eta)$.  To emphasize the dependence of this module on $\eta$, 
one could write it as $(\eta A)^\perp/A$, where $(\eta A)^\perp 
\subset H^0_*(\mathcal O_Z)$ means the orthogonal complement of $\eta 
A \subset H^0_*(\omega_Z)$ with respect to the Serre duality pairing. 
Now $(\eta A)^\perp /A$ is dual to $(\omega_A/\eta A)$, and it is also 
self-dual with a shift.  Consequently if $\eta \in \omega_A$ is of 
degree $t$, then the corresponding deficiency module is 
\[ 
(\eta A)^\perp / A \cong (\omega_A/\eta A)' \cong (\omega_A/\eta 
A)(t). 
\] 
 
In Theorem \ref{pfaff.points} we split the deficiency module in half, 
and put a Lagrangian subhalf in $\mathcal F$ and the quotient half in 
$\mathcal F^*(t-4)$.  The Pfaffian resolutions of $\mathcal O_Z$ are 
thus classified up to symmetric homotopy equivalence by pairs $(\eta, 
L/A)$ with $\eta \in \omega_A$ generating the sheaf $\omega_Z$, and 
with $L/A \subset (\eta A)^\perp /A$ a Lagrangian submodule. 
\end{remark} 
 
An alternative strategy for dealing with this deficiency module is to 
construct a diagram of the form of \eqref{dual.diag} in Theorem 
\ref{first.v} (we write $\mathcal O := \mathcal O_{\mathbb 
P^3}$ to try to stay inside the margins): 
\begin{small}%
\begin{equation}
\label{sublagr.dual} 
\begin{diagram}
0 & \rTo & \mathcal O(t-4) & \rTo & \mathcal G & \rTo ^\psi & 
\bigoplus \mathcal O(-a_i) & \rTo & \mathcal O & \rOnto & \mathcal 
O_Z \\ 
&& \dSame && \dTo <\phi && \dTo >{\phi^*} && \dSame && \dTo >\eta 
<\cong \\ 
0 & \rTo & \mathcal O(t-4) & \rTo & \bigoplus \mathcal O(a_i+t-4) & 
\rTo _{-\psi^*} & \mathcal G^*(t-4) & \rTo & \mathcal O & \rOnto & 
\omega_Z(t) 
\end{diagram} 
\end{equation}%
\end{small}
with $\bigoplus \mathcal O(-a_i)$ corresponding to a minimal set of 
generators of the homogeneous ideal of $Z$, with $H^2_*(\mathcal G) 
\cong (\eta A)^\perp / A$, the deficiency module, and with 
$H^1_*(\mathcal G) = 0$.

We now give examples both of Pfaffian resolutions as in 
\eqref{pfaff.resol} which split the deficiency module, and of 
resolutions as in \eqref{sublagr.dual} in the form of Theorem 
\ref{first.v} which gather the deficiency module up in one piece. 
 
\subsubsection*{Example 1: one point} 
 
Consider a single rational point $Q$.  Its geometry is simple, but we
can make its algebra complicated.
 
The canonical module of $Q$ is $\omega_A \cong \bigoplus _{n\geq 1} 
H^0(\omega_Q(n))$.  If we pick a nonzero $\eta$ of degree $1$, then it 
generates $\omega_A$, and its deficiency module vanishes.  The 
constructions described above both lead unsurprisingly to the Koszul 
resolution 
\[ 
0 \to \mathcal O_{\mathbb P^3}(-3) \to \mathcal O_{\mathbb 
P^3}(-2)^{\oplus 3} \to \mathcal O_{\mathbb P^3}(-1)^{\oplus 3} \to 
\mathcal O_{\mathbb P^3} \to \mathcal O_Q \to 0. 
\] 
 
However, if we let $\eta \in \omega_A$ be a nonzero element of degree 
$2$, then the deficiency module is $k$ concentrated in degree $-1$, 
and the construction \eqref{sublagr.dual} yields a diagram 
\begin{diagram} 
0 & \rTo & \mathcal O_{\mathbb P^3}(-2) & \rTo & \Omega^2_{\mathbb 
P^3}(1) & \rTo & \mathcal O_{\mathbb P^3}(-1)^{\oplus 3} & \rTo & 
\mathcal O_{\mathbb P^3} & \rOnto & \mathcal O_Q \\ 
&& \dSame && \dTo && \dTo && \dSame && \dTo >\eta <\cong \\ 
0 & \rTo & \mathcal O_{\mathbb P^3}(-2) & \rTo & \mathcal O_{\mathbb 
P^3}(-1)^{\oplus 3} & \rTo & \Omega_{\mathbb P^3}(1) & \rTo & \mathcal 
O_{\mathbb P^3} & \rOnto & \omega_Q(2). 
\end{diagram} 

More generally, if we let $\eta \in \omega_A$ be a nonzero element of 
degree $t$, then the deficiency module is $\bigoplus_{n=-(t-1)}^{-1} 
H^0(\mathcal O_Q(n))$, and the construction \eqref{sublagr.dual} 
yields 
\begin{diagram} 
0 & \rTo & \mathcal O_{\mathbb P^3}(t-4) & \rTo & \mathcal F_t^*(t-4) 
& \rTo & \mathcal O_{\mathbb P^3}(-1)^{\oplus 3} & \rTo & \mathcal 
O_{\mathbb P^3} & \rOnto & \mathcal O_Q \\ 
&& \dSame && \dTo && \dTo && \dSame && \dSame \\  
0 & \rTo & \mathcal O_{\mathbb P^3}(t-4) & \rTo & \mathcal O_{\mathbb 
P^3}(t-3)^{\oplus 3} & \rTo & \mathcal F_t & \rTo & \mathcal 
O_{\mathbb P^3} & \rOnto & \mathcal O_Q 
\end{diagram} 
with $\mathcal F_t$ a rank $3$ locally free sheaf which is the 
sheafification of the kernel of the presentation of the deficiency 
module: 
\[ 
0 \to \mathcal F_t \to \mathcal O_{\mathbb P^3} \oplus \mathcal 
O_{\mathbb P^3}(t-2)^{\oplus 3} \to \mathcal O_{\mathbb P^3}(t-1) \to 
0. 
\] 
 
If one lets $\eta \in \omega_A$ be a nonzero element of degree $3$, 
then applying the methods of Theorem \ref{pfaff.points} yields a 
resolution which one recognizes as the Koszul complex associated to 
the zero locus of a section of the rank $3$ bundle $\mathcal 
T_{\mathbb P^3}(-1)$: 
\[ 
0 \to \mathcal O_{\mathbb P^3}(-1) \to \Omega^2_{\mathbb P^3}(2) \to 
\Omega_{\mathbb P^3}(1) \to \mathcal O_{\mathbb P^3} \to \mathcal O_Q 
\to 0. 
\]

\subsubsection*{Example 2: three points} 
 
If $Z$ is the union of three non-collinear rational points, then the 
module $\omega_A$ has two generators of degree $0$, and $Z$ is not 
arithmetically Gorenstein.  If we pick a general $\eta \in \omega_A$ 
of degree $0$, then the deficiency module is $k$, concentrated in 
degree $0$, and the construction \eqref{sublagr.dual} yields a diagram 
(in which we again write $\mathcal O := \mathcal O_{\mathbb P^3}$ in 
order to simplify the notation): 
\begin{small} 
\begin{diagram} 
0 & \rTo & \mathcal O(-4) & \rTo & \mathcal O(-3) \oplus 
\Omega^2_{\mathbb P^3} & \rTo & \mathcal O(-2)^{\oplus 3} \oplus 
\mathcal O(-1) & \rTo & \mathcal O & \rOnto & \mathcal O_Z \\ 
&& \dSame && \dTo && \dTo && \dSame && \dSame \\ 
0 & \rTo & \mathcal O(-4) & \rTo & \mathcal O(-3) \oplus \mathcal 
O(-2)^{\oplus 3} & \rTo & \Omega_{\mathbb P^3} \oplus \mathcal O(-1) & 
\rTo & \mathcal O & \rOnto & \mathcal O_Z 
\end{diagram} 
\end{small} 

If we pick a general $\eta \in \omega_A$ of degree $1$, then the 
deficiency module $(\omega_A/\eta A)(1)$ is of length $4$, 
concentrated in degrees $0$ and $-1$, and the methods of Theorem 
\ref{pfaff.points} yield a symmetric resolution (with alternating 
middle map $\psi$): 
\[ 
0 \to \mathcal O(-3) \to \mathcal \Omega^2_{\mathbb P^3}(1)^{\oplus 2} 
\oplus \mathcal O(-2) \xrightarrow{\,\psi\,} \Omega_{\mathbb P^3}^{\oplus 
2} \oplus \mathcal O(-1) \to \mathcal O \to \mathcal O_Z \to 0. 
\]

\section{Some weakly subcanonical subschemes} 
\label{examples} 
 
In this section we give some examples of weakly subcanonical
subschemes.  These are examples of subschemes $Z \subset X$ which
satisfy conditions \cond1-\cond2 of the definition of a strongly
subcanonical subscheme but fail one or both of conditions
\cond3-\cond4.  Thus the Serre construction (in codimension $2$) and
our Theorem \ref{converse} (in codimension $3$) fail for these
subschemes.

\subsection{A weakly subcanonical curve} 
 
We construct a subcanonical curve $C \subset \mathbb P^1 \times 
\mathbb P^n$ for $n \geq 2$ which fails the lifting condition \cond4 
of Definition \ref{strong.def}.

Let $C$ be a nonsingular projective curve of genus $2$ over an 
algebraically closed field $k$, let $P$ be one of its Weierstrass 
points, and let $D$ be a divisor of degree $4$ on $C$.  A 
base-point-free pencil in the linear system of divisors $\linsys D$ 
defines a map $f : C \to \mathbb P^1$, and a base-point-free net in 
$\linsys {D+P}$ defines a map $g : C \to \mathbb P^2$.  Composing $g$ 
with a linear embedding $\mathbb P^2 \hookrightarrow \mathbb P^n$ 
gives a map $h : C \to \mathbb P^n$.  Let $i := (f,h): C \to \mathbb 
P^1 \times \mathbb P^n$.  If the linear systems are chosen 
sufficiently generally, then $i$ is an embedding.

The restriction to $C$ of a line bundle $\mathcal O_{\mathbb P^1
\times \mathbb P^n} (a,b)$ is $\mathcal O_C((a+b)D+bP)$.  So the
canonical bundle $\omega_C \cong \mathcal O_C(2P)$ is the restriction
of $\mathcal O_{\mathbb P^1 \times \mathbb P^n} (-2,2)$.  If the class
of $D-4P$ in $\Pic^0(C)$ is not torsion, then $\mathcal O_{\mathbb P^1
\times \mathbb P^n} (-2,2)$ is the only line bundle on $\mathbb P^1
\times \mathbb P^n$ whose restriction is $\omega_C$.  Hence the
subcanonical curve $C \subset \mathbb P^1 \times \mathbb P^n$ will
definitely fail such structure theorems as the Serre construction or
Theorem \ref{converse} if the lifting condition \cond4 of Definition
\ref{strong.def} fails for the isomorphism $\eta : \omega_C \cong
\mathcal O_{\mathbb P^1 \times \mathbb P^n} (-2,2) \rest C$.
 
By \eqref{cond.dual} this failure is equivalent to the nonvanishing of 
the composite map 
\begin{equation} 
\label{dual.lifting} 
H^1(\mathbb P^1 \times \mathbb P^n, \mathcal O(-2,2)) 
\xrightarrow{\text{\rm rest}} H^1(C,\mathcal O(-2,2) \rest C) 
\xrightarrow[\cong]{\eta} H^1(C,\omega_C) \xrightarrow[\cong]{\tr} k. 
\end{equation} 
Now the image of $g : C \to \mathbb P^2$ is a singular quintic plane 
curve.  If we resolve the singularities, then $g$ factors as an 
embedding followed by the blowdown $C \hookrightarrow \widetilde 
{\mathbb P}^2 \to \mathbb P^2$.  The composite map of 
\eqref{dual.lifting} now factors through the diagram 
\begin{footnotesize} 
\begin{diagram} 
H^1(\mathbb P^1 \times \mathbb P^n, \mathcal O(-2,2)) & \rOnto 
^{\qquad} & H^1(\mathbb P^1 \times \mathbb P^2, \mathcal O(-2,2)) & 
\lTo ^ \cong & H^1(\mathbb P^1 \times \widetilde{\mathbb P}^2, 
\mathcal O(-2,2)) \cong k^6 \\ 
&& \dTo <\alpha && \dTo >\beta  \\ 
&& k \cong H^1(C, \omega_C) & \lTo ^\gamma & H^1(\mathbb P^1 \times C, 
\mathcal O(-2,2(D+P))) \cong k^9 
\end{diagram} 
\end{footnotesize}%
The lifting condition \cond4 fails if and only if $\alpha$ is 
surjective, hence if and only if $\im(\beta) \not\subset 
\ker(\gamma)$. 
 
Now $\gamma$ is part of the long exact sequence of cohomology for 
\[ 
0 \to \mathcal O_{\mathbb P^1 \times C}(-3,D+2P) \to \mathcal 
O_{\mathbb P^1 \times C}(-2,2D+2P) \to \omega_C \to 0. 
\] 
So $\gamma$ is surjective, and $\ker(\gamma) \subset k^9$ is a 
hyperplane.   
 
The map $g : C \to \mathbb P^2$ is defined using a three-dimensional 
subspace $U_3 \subset V_4 := H^0(C,\mathcal O_C(D+P))$.  The complete 
linear system embeds $C \hookrightarrow \mathbb P^3$ as a curve of 
degree $5$ and genus $2$ contained in a unique quadric surface $Q$. 
Then $\beta$ is the natural map from $S^2 U_3 \cong k^6$ to $S^2 V_4 / 
\langle Q \rangle \cong k^9$.  Now we have a range of choices for the 
subspace $U_3 \subset V_4$ which vary in a Zariski open subset of 
$\mathbb P^3 = \mathbb P(V_4^*)$.  Hence we have a family of possible 
subspaces $S^2 U_3 \subset S^2 V_4$ whose different members are not 
all contained in any fixed hyperplane of $S^2 V_4$.  So if we choose a 
general $U_3 \subset V_4$, then $S^2 U_3 = \im(\beta)$ is not 
contained in the hyperplane $\ker(\gamma) \subset S^2 V_4 / \langle Q 
\rangle$.  In that case, $C \subset \mathbb P^1 \times \mathbb P^n$ is 
a subcanonical curve which fails the lifting condition \cond4.

\subsection{Singular points} 
 
Examples can easily be given of subcanonical subschemes $Z \subset X$
which are not covered by our construction because the finite
projective dimension condition \cond3 of Definition \ref{strong.def}
breaks down.  This may happen at the same time that the lifting
condition \cond4 breaks down, or it may happen independently.  If
\cond4 holds but \cond3 breaks down, $\mathcal O_Z$ will still have
resolutions fitting into diagrams such as \eqref{dual.diag} of Theorem
\ref{first.v}, except that $\mathcal E$ or $\mathcal F$ will not be
locally free.
 
For instance if $X \subset \mathbb P^4$ is a singular hypersurface of 
degree $d$, and $P \in X$ is a singular point, then $P$ is indeed 
subcanonical, but condition \cond3 fails because $\mathcal O_P$ is of 
infinite local projective dimension over $\mathcal O_X$.  There exist 
isomorphisms $\eta : \mathcal O_P \cong \shExt^3_{\mathcal 
  O_{X}}(\mathcal O_P,\mathcal O_{X}(\ell))$ for all $\ell \in \mathbb 
Z$, but these satisfy condition \cond4 if and only if $\ell \geq d-4$. 
 
Similarly, if $D$ is a line in $\mathbb P^5$, and $Y \subset \mathbb 
P^5$ a hypersurface containing $D$ which is singular in at least one 
point of $D$, then condition \cond3 fails for $D \subset Y$, but all 
the other conditions hold (since $H^3(Y, \omega_{Y}(2)) = 0$).  So 
although $D \subset Y$ may be obtained as a degeneracy locus of a pair 
of Lagrangian subsheaves of a twisted orthogonal bundle on $Y$, at 
least one of the Lagrangian subsheaves is not locally free.

\subsection{A nonseparated example} 
 
We now give an example where there is no real choice about the 
$\eta$ (because $H^0(\omega_Z) = k$ and there are no twists), 
where conditions \cond1-\cond3 hold, but where condition \cond4 
fails.  The real reason for the failure in this example is that 
we are doing something silly on a nonseparated scheme.  But the 
interesting thing is that the cohomological obstruction \cond4 
is able to detect our misbehavior. 
 
Let $X$ be the nonseparated scheme consisting of two copies $\mathbb
A^3$ glued together along $\mathbb A^3 - \{0\}$.  In other words, $X$
is $\mathbb A^3$ with the origin doubled up.  Let $P' \in X$ be one of
the two origins.  It is a subcanonical subscheme of $X$ of codimension
$3$ of finite local projective dimension, i.e.\ it satisfies
conditions \cond1-\cond3 of Definition \ref{strong.def}.  We claim
that it does not satisfy condition \cond4.
 
The problem is to compute the map 
\begin{equation} 
\label{ext3} 
\Ext^3_{\mathcal O_X}(\mathcal O_{P'},\mathcal O_X) \to 
\Ext^3_{\mathcal O_X}(\mathcal O_X,\mathcal O_X) = H^3(X,\mathcal 
O_X). 
\end{equation} 
We will use the following notation: $U',U'' \subset X$ are the two 
copies of $\mathbb A^3$; for $\alpha=1,2,3$ let $U_\alpha \subset X$ 
be the open locus where $x_\alpha \neq 0$; and let $U_{\alpha\beta} := 
U_\alpha \cap U_\beta$, and $U_{123} := U_1 \cap U_2 \cap U_3$.  For 
any inclusion of an affine open subscheme $U \subset X$, we denote by 
$i_! \mathcal O_U$ the extension by zero of $\mathcal O_U$ to all of 
$X$.  We will use the same letter $i_!$ whatever the $U$. 
 
Then $\mathcal O_X$ and $\mathcal O_{P'}$ have resolutions of the form 
\begin{diagram} 
0 & \rTo & i_! \mathcal O_{U_{123}} & \rTo & \bigoplus_{\alpha < 
\beta} i_! \mathcal O_{U_{\alpha\beta}} & \rTo & \bigoplus_\alpha i_! 
\mathcal O_{U_\alpha} & \rTo & i_! \mathcal O_{U'} \oplus i_! \mathcal 
O_{U''} & \rOnto & \mathcal O_X \\  
&& \dTo && \dTo && \dTo && \dTo && \dTo \\  
0 & \rTo & i_! \mathcal O_{U'} & \rTo & i_! \mathcal O_{U'}^{\oplus 
3}& \rTo & i_!  \mathcal O_{U'}^{\oplus 3} & \rTo & i_!  \mathcal 
O_{U'} & \rOnto & \mathcal O_{P'} 
\end{diagram} 
The horizontal maps in the first row are more or less taken from a \v 
Cech resolution, while those from the second row are from a Koszul 
resolution.  The vertical maps are, from left to right, 
\begin{align*} 
& \begin{pmatrix} \frac 1 {x_1 x_2 x_3} \end{pmatrix},  
&& \begin{pmatrix} \frac 1 {x_1 x_2} & 0 & 0 \\ 0 & \frac 1 {x_1 x_3} 
& 0 \\ 0 & 0 & \frac 1 {x_2 x_3} \end{pmatrix}, 
&& \begin{pmatrix} \frac 1 {x_1} & 0 & 0 \\ 0 & \frac 1 {x_2} & 0 \\ 0 
& 0 & \frac 1 {x_3} \end{pmatrix}, 
&& \begin{pmatrix} 1 & 0 \end{pmatrix}. 
\end{align*} 
If we apply $\Hom_{\mathcal O_X}({-},\mathcal O_X)$ to the 
resolutions, we get complexes which compute the $\Ext^p_{\mathcal 
O_X}(\mathcal O_{P'},\mathcal O_X)$ and the $H^p(X,\mathcal O_X)$. 
(This is because $\mathcal O_X$ is quasi-coherent, and the $i_! 
\mathcal O_U$ are extensions by zero of locally free sheaves on affine 
open subschemes.)  Writing $R := k[x_1,x_2,x_3]$, these complexes are 
\begin{diagram} 
0 & \rTo & R & \rTo & R^{\oplus 3} & \rTo & R^{\oplus 3} & \rTo & R & 
\rTo & 0 \\ 
&& \dTo && \dTo && \dTo && \dTo \\ 
0 & \rTo & R \oplus R & \rTo & \bigoplus_{\alpha} R[x_\alpha^{-1}] & 
\rTo & \bigoplus_{\alpha < \beta} R[x_\alpha ^{-1}x_\beta^{-1}] & \rTo 
& R[x_1^{-1} x_2^{-1} x_3^{-1}] & \rTo & 0 
\end{diagram} 
with Koszul and \v Cech horizontal arrows.  The vertical arrows are as 
before, but transposed and in the reverse order.   
 
The map \eqref{ext3} which we wish to compute may now be identified as
$k \hookrightarrow H^3_{\mathfrak m}(R)$.  This map sends a nonzero
$\eta$ to a nonzero multiple of socle element $x_1^{-1} x_2^{-1}
x_3^{-1}$ of $H^3_{\mathfrak m}(R)$.  So $P' \subset X$ fails
condition \cond4 of the definition of a strongly subcanonical
subscheme.
 
Actually, this example has additional pathologies which prevent 
$\mathcal O_{P'}$ from having a locally free resolution, 
whether symmetric or otherwise.  For $\mathcal O_{P'}$ does not 
even have a locally free presentation.  Indeed, for reasons of 
depth, any map $\mathcal E \to \mathcal F$ between locally free 
sheaves on $X$ is determined by what happens outside the two 
origins.  So the cokernel of such a map has the same fiber at 
the two origins.  Therefore the cokernel cannot be $\mathcal 
O_{P'}$.

\section{Codimension one sheaves} 
\label{codimension.one} 
 
In this section we consider the analogues of the results in the 
previous sections for (skew)-symmetric sheaves of codimension $1$.  We 
include necessary and sufficient conditions for such sheaves on 
$\mathbb P^N$ to have locally free resolutions which are genuinely 
(skew)-symmetric, similar to those in Walter \cite{Walter}.  We also 
prove that any such sheaf on a quasi-projective variety has a 
resolution which is (skew)-symmetric up to quasi-isomorphism, in 
analogy with Theorem \ref{converse}.  We finish the section with 
several examples.

In this section we will suppose that the characteristic is 
not $2$, although all the theorems have variants which are valid 
in characteristic $2$. 
 
\subsection{Symmetric sheaves of codimension 1} 
 
Suppose $\mathcal F$ is a coherent sheaf on a scheme $X$ which is of 
finite local projective dimension and perfect of codimension $1$. 
This means that locally $\mathcal F$ has free resolutions $0 \to 
\mathcal L_1 \to \mathcal L_0 \to \mathcal F \to 0$ such that the dual 
complex $0 \to \mathcal L_0^* \to \mathcal L_1^* \to 
\shExt^1_{\mathcal O_X}(\mathcal F,\mathcal O_X) \to 0$ is also exact. 
The operation 
\[ 
\mathcal F \rightsquigarrow \mathcal F^\vee := \shExt^1_{\mathcal O_X} 
(\mathcal F, \mathcal O_X) 
\] 
provides a duality on the category of such sheaves.  A {\em symmetric 
sheaf of codimension $1$} is a pair $(\mathcal F, \alpha)$ where 
$\mathcal F$ is a sheaf as above, and $\alpha : \mathcal F \to 
\mathcal F^\vee(L)$ is an isomorphism which is symmetric in the sense 
that $\alpha = \alpha^\vee$.  (Here $L$ is some line bundle on $X$.) 
{\em Skew-symmetric sheaves of codimension $1$} on $X$ are defined 
similarly.

\subsection{Symmetric resolutions in codimension $1$} 
 
Resolutions of codimension $1$ symmetric sheaves on $\mathbb 
P^3$ have been studied fairly extensively by Barth 
\cite{Barth}, Casnati-Catanese \cite{CC}, and Catanese 
\cite{Catanese2} \cite{Catanese}, in the context of surfaces 
with even sets of nodes and by Kleiman-Ulrich \cite{KU} in the 
context of self-linked curves.  The next theorem, conjectured 
by Barth and Catanese, was proven by Casnati-Catanese for 
symmetric sheaves on $\mathbb P^3$ (\cite{CC} Theorem 0.3). 
They also remarked (\cite{CC} Remark 2.2) that essentially the 
same proof works for codimension $1$ symmetric sheaves on any 
$\mathbb P^n$, which is true as long as one remembers to 
include in one's statement a parity condition analogous to that 
in Walter \cite{Walter} Theorem 0.1.  For a case where the 
parity condition fails, see Example \ref{threefold} below.

\begin{theorem}[\cite{CC} \cite{Catanese} with correction] 
\label{CCatanese} 
Let $k$ be an algebraically closed field of characteristic different 
from $2$.  Suppose that $(\mathcal F,\alpha)$ is a symmetric sheaf of 
codimension $1$ on $\mathbb P^{n}_k$, with $\alpha : \mathcal F 
\xrightarrow{\sim} \mathcal F^\vee(\ell-n-1)$.  Then $\mathcal F$ has 
a symmetric resolution, i.e.\ a locally free resolution of the form 
\[ 
0 \to \mathcal G \xrightarrow{\,f\,} \mathcal G^*(\ell-n-1) \to 
\mathcal F \to 0 
\] 
with $f$ symmetric, if and only if the following parity condition 
holds\textup{:} if $n\equiv 1 \pmod{4}$ and $\ell$ is even, then $\chi 
(\mathcal F (-\ell/2))$ is also even. 
\end{theorem} 
 
A higher-codimension generalization of this theorem is proven in 
our paper \cite{EPW3}. 
 
As the parity condition indicates, symmetric sheaves do not always 
possess symmetric resolutions.  The following structure theorem, 
analogous to Theorem \ref{main.result}, shows that they do still have 
locally free resolutions which are symmetric up to quasi-isomorphism.

\begin{theorem}
\label{quasisym}
Let $X$ be a quasiprojective scheme over a Noetherian ring and
let $\mathcal F$ be a coherent sheaf.
The following are equivalent:

\textup{(a)} $\mathcal F$ is perfect of codimension one, and 
there exists a line bundle $L$ on $X$ and an
isomorphism $\alpha:\mathcal F\to \mathcal F^\vee(L)$ making
$(\mathcal F,\alpha)$ a symmetric sheaf.

\textup{(b)} The sheaf $\mathcal F$ has symmetrically quasi-isomorphic
locally free resolutions
\begin{equation} 
\label{quasisym.diag} 
\begin{diagram}
0 & \rTo & \mathcal G  & \rTo ^\psi & \mathcal H & \rTo & \mathcal F 
& \rTo & 0 \\  
&& \dTo <\phi && \dTo >{\phi^*} && \dTo >\alpha 
<\cong \\ 
0 & \rTo & \mathcal H^*(L) & \rTo ^{\psi^*} & \mathcal 
G^*(L)& \rTo & \mathcal F^\vee(L) & \rTo & 0 
\end{diagram} 
\end{equation}
%
with $\phi^*\psi : \mathcal G \to \mathcal G^*(L)$ a symmetric 
map and $L$ a line bundle on $X$.

\textup{(c)} There exists a line bundle $L$ on $X$ and a Lagrangian subbundle
of a twisted symplectic bundle 
\[
\mathcal G\rInto^{\sm {\psi\\ \phi}}\mathcal H \oplus \mathcal H^*(L)
\] 
such that 
\begin{diagram}
0 & \rTo & \mathcal G  & \rTo ^\psi & \mathcal H & \rTo & \mathcal F 
& \rTo & 0 \\ 
\end{diagram}
is a resolution of $\mathcal F$ fitting into a commutative diagram
as in \eqref{quasisym.diag}.
\end{theorem}

\begin{proof} The only delicate part is (a) $\Rightarrow$ (b).
Because $\mathcal F$ is locally Cohen-Macaulay of codimension $1$, it 
has a locally free resolution $0 \to \mathcal P_1 \to \mathcal P_0 \to 
\mathcal F \to 0$.  The symmetric isomorphism $\alpha$ corresponds to 
a morphism in the derived category 
\begin{diagram} 
S^2(\mathcal P_\cpx) :  & \qquad &  
0 & \rTo & \Lambda^2 \mathcal P_1 & \rTo & \mathcal P_1 \otimes 
\mathcal P_0 & \rTo & S^2 \mathcal P_0 & \rTo & 0 \\ 
\dDotsto >\alpha && && \dDotsto && \dDotsto && \dDotsto  \\ 
L[1]: && 0 & \rTo & 0 & \rTo & L & \rTo & 0 & \rTo & 0  
\end{diagram} 
This $\alpha$ is a member of the hyperext $\hExt^1_{\mathcal O_X} 
(S^2(\mathcal P_\cpx), L)$, which in turn is the abutment of the 
hyperext spectral sequence 
\[ 
E_1^{pq} = \Ext^q_{\mathcal O_X}((S^2(\mathcal P_\cpx))_{p} , L) 
\Longrightarrow \hExt^{p+q}_{\mathcal O_X} (S^2(\mathcal P_\cpx), L). 
\] 
The differentials $d_1$ define complexes (indexed by $p=0,1,2$) 
\[ 
0 \to \Ext^q_{\mathcal O_X}(S^2 \mathcal P_0,L) \xrightarrow{d_1} 
\Ext^q_{\mathcal O_X} (\mathcal P_1 \otimes 
\mathcal P_0,L) \xrightarrow{d_1} \Ext^q_{\mathcal O_X}(\Lambda^2 
\mathcal P_1,L) \to 0 
\] 
whose cohomology groups are the $E_2^{pq}$.  In particular, $E_2^{10}$ 
is the space of homotopy classes of chain maps $S^2(\mathcal P_\cpx) 
\to L[1]$.  Hence $\alpha$ will be the class of an honest chain map if 
and only if it comes from $E_2^{10}$.  However, according to the 
$5$-term exact sequence 
\[ 
0 \to E_2^{10} \to \hExt^1_{\mathcal O_X} (S^2(\mathcal P_\cpx), L) 
\to E_2^{01} \to \dotsb, 
\] 
the obstruction lies in $E_2^{01} \subset \Ext^1_{\mathcal O_X} (S^2 
\mathcal P_0, L)$.  As in the proof of Theorem \ref{converse}, this 
obstruction may be nonzero, but it can be killed by pulling back along 
a suitable epimorphism $\mathcal H \twoheadrightarrow \mathcal P_0$ 
(cf.\ Lemma \ref{quasiproj}).  The proof may now be completed with 
arguments taken from the proof of Theorem \ref{converse}. 
\end{proof} 
 
We now use this theorem to construct an example of a symmetric 
codimension $1$ sheaf on $\mathbb P^5$ for which the parity condition 
of Theorem \ref{CCatanese} fails.   The construction is similar to the 
main examples of \cite{EPW2}.

\begin{example} 
\label{threefold} 
Let $V = H^0(\mathcal O_{\mathbb P^5}(1))^*$.  The exterior product 
defines a symplectic form on the $20$-dimensional vector space 
$\Lambda^3 V$, which makes the trivial bundle $\Lambda^3 V \otimes 
\mathcal O_{\mathbb P^5}$ into a symplectic bundle.  If $W, W^*\subset 
\Lambda^3 V$ are general Lagrangian subspaces, then we can identify 
the symplectic vector space $\Lambda^3 V$ with hyperbolic symplectic 
vector space $W\oplus W^*$.  Moreover, one can see that 
$\Omega^3_{\mathbb P^5}(3)$ is a Lagrangian subbundle of $\Lambda^3 V 
\otimes \mathcal O_{\mathbb P^5}$ (cf.\ \cite{EPW2}, \S 5). 
The construction of Theorem \ref{quasisym} then produces a symmetric 
codimension $1$ sheaf $\mathcal F$ on $\mathbb P^5$ with resolutions 
\begin{equation}
\label{example.diag} 
\begin{diagram}
0 & \rTo & \Omega^3_{\mathbb P^5}(3) & \rTo^\psi  
& W \otimes \mathcal O_{\mathbb P^5} &  
\rTo & \mathcal F & \rTo & 0 \\ 
&& \dTo <\phi && \dTo >{\phi^*} && \dTo >\alpha <\cong \\ 
0 & \rTo & W^*\otimes \mathcal O_{\mathbb P^5} &  
\rTo ^{\psi^*} &\Omega^2_{\mathbb P^5}(3) & \rTo & 
\mathcal F^\vee & \rTo & 0. 
\end{diagram}%
\end{equation}%
The sheaf $\mathcal F$ fails the parity condition of Theorem 
\ref{CCatanese} because $\ell=6$ and $\chi(\mathcal F(-3))=1$. 
 
The geometry of the sheaf $\mathcal F$ is best explained using the 
degeneracy loci of the Lagrangian subbundles $\Omega^3_{\mathbb 
P^5}(3)$ and $W^* \otimes \mathcal O_{\mathbb P^5}$ of $\Lambda^3 V 
\otimes \mathcal O_{\mathbb P^5}$: 
\[ 
D_i := \{ x \in \mathbb P^5 \mid \dim \left[ \Omega^3_{\mathbb 
P^5}(3)(x) \cap W^* \right] \geq i \}. 
\] 
The sheaf $\mathcal F$ is supported on the sextic fourfold $D_1$.  If 
$W^*$ is general, then $D_1$ is smooth (cf.\ \cite{EPW2}, Theorem 2.1 
and the discussion following it)  except 
along the surface $D_2$ where it has $A_1$ singularities with local 
equations $x_1^2 + x_2^2 + x_3^2 = 0$.  The surface $D_2$ is of degree 
$40$ according to the formulas of Fulton-Pragacz \cite{FP} (6.7). 
 
Now choose a general $9$-dimensional subspace $U$ of the 
$10$-dimensional space $W^*$.  Then the composite map $U \otimes 
\mathcal O_{\mathbb P^5} \hookrightarrow \Lambda^3 V \otimes \mathcal 
O_{\mathbb P^5} \twoheadrightarrow \Omega_{\mathbb P^5}^2 (3)$ 
degenerates in codimension $2$ along a threefold $Y$ of degree $18$. 
Since 
\[ 
Y = \{ x \in \mathbb P^5 \mid \dim \left[ \Omega^3_{\mathbb P^5}(3)(x) 
\cap U \right] \geq 1 \}, 
\] 
we have $D_2 \subset Y \subset D_1$.  Moreover, $\mathcal F \cong 
\mathcal I_{Y/D_1}(6)$.  In addition, $Y$ is self-linked by the complete 
intersection of $D_1$ and of another sextic hypersurface corresponding 
to another Lagrangian subspace of $\Lambda^3 V$ containing $U$. 
\end{example} 
 
\subsection{Skew-symmetric sheaves of codimension one} 
 
Analogues of Theorems \ref{CCatanese} and \ref{quasisym} 
hold for skew-symmetric sheaves of codimension $1$.  The 
only significant change is in the parity condition of Theorem 
\ref{CCatanese}, which in the skew-symmetric case has the form: ``if 
$n\equiv 3 \pmod{4}$ and $\ell$ is even, then $\chi (\mathcal F 
(-\ell/2))$ is also even.''  We leave the exact formulation of these 
results to the reader. 
 
If $S \subset \mathbb P^3$ is a smooth surface of degree $d$, then its 
cotangent bundle $\Omega_S$ is a skew-symmetric sheaf of codimension 
$1$ on $\mathbb P^3$ with twist $\ell = 0$.  Since $\chi(\Omega_S) = 
-h^{11}(S) 
\equiv d \pmod 2$, this skew-symmetric sheaf fails the parity 
condition when $d$ is odd.

\end{document}